\documentclass[12pt]{article}
\usepackage{stmaryrd}
\usepackage{amssymb}
\usepackage{color, amsmath,amssymb, amsfonts, amstext,amsthm, latexsym}

\usepackage{amssymb, epsfig, amssymb, latexsym}
\usepackage{amsmath}
\usepackage{graphicx}
\usepackage{longtable}
\usepackage{graphicx}
\usepackage{subfigure}
\usepackage{appendix}
\usepackage[symbol]{footmisc}
\usepackage[colorlinks=false, linktocpage=true]{hyperref}
\usepackage[noblocks]{authblk}
\numberwithin{equation}{section} \topmargin -0.4in

\renewcommand{\phi}{\varphi}

\begin{document}
\title{L\'evy noise induced escape in the Morris-Lecar model}

\author[1,2]{Yancai Liu}
\affil[1]{{School of Mathematics and Statistics, Huazhong University of Science and Technology,
 Wuhan 430074, China}}
\affil[2]{{Center for Mathematical Sciences, Huazhong University of Science and Technology, Wuhan 430074, China}}
\author[3 \footnote{Corresponding author: e-mail: cairui5876@163.com.}]{Rui Cai}
\affil[3]{{School of Science, Hubei University of Technology, Wuhan 430068, China}}
\author[4]{Jinqiao Duan}
\affil[4]{{Department of Applied Mathematics, Illinois Institute of Technology,
 Chicago, IL 60616, USA}}

\date{\today}
\maketitle

\begin{abstract}

The phenomenon of an excitable system producing
a pulse  under external or internal stimulation may be interpreted as a stochastic escape problem. This work addresses this issue by examining the Morris-Lecar neural model driven by symmetric $\alpha$-stable L\'evy motion (non-Gaussian noise) as well as Brownian motion (Gaussian noise). Two deterministic quantities: the first escape probability  and the mean first exit time, are adopted to analyse the state transition from the resting state to the excited state and  the stability of this stochastic  model.
Additionally, a recent geometric concept, the  stochastic basin of attraction  is used to explore the basin stability of the escape region.
Our main results include: (i) the larger L\'evy motion index  with smaller jump magnitude and the relatively small noise intensity are conducive for the Morris-Lecar model to produce pulses; (ii) a smaller noise intensity  and a larger L\'evy motion index  make the mean first exit time longer, which means the stability of the resting state can be enhanced in this case; (iii) the effect of ion channel noise is more pronounced on the stochastic Morris-Lecar model than the current noise.
 This work provides some numerical simulations about the impact of non-Gaussian, heavy-tailed, burst-like fluctuations on excitable systems such as the Morris-Lecar system.

\end{abstract}


\section{Introduction}
In recent years, there has been an increasing interest in the effects of noise in neuroscience. Because of the noisy environment that neurons live in, there are many sources of noise in the neuronal systems.  These noisy sources include, for instance,
random attacks caused by spontaneous release of neurotransmitters, synaptic noise from spontaneous postsynaptic potentials, small fluctuations in the electrical potential across the nerve-cell membrane,  and the opening and closing of ion channels. Noise may induce various phenomena, such as oscillations \cite{Justus2010,wang2017}, chaos-like behaviors  \cite{gao1999can}, state transitions \cite{basu2001spontaneously}, stochastic resonance \cite{wang2016, mcdonnell2008stochastic,Dybiec2009levy}, and spatial coherence resonance \cite{Tao2017resonance,Gu2013b,Perc2005spatial}.

 In this paper,  we study the Morris-Lecar (ML) neuron model under the disturbance of (non-Gaussian) L\'evy noise as well as (Gaussian) Brownian  noise. As a simplified version of the Hodgkin-Huxley (HH) system, the ML model was first introduced to account for the electrical activities of the giant barnacle muscle fibers in invertebrates \cite{ML1981}. Since then, it becomes a canonical neuronal model because it can display two different forms of neuronal excitability behaviors under various
parameter regimes. One form of excitability  is not sensitive to external stimulation intensity, the discharge starting frequency can be very low, and the discharge range is relatively wide. This is  called type I excitability. The type II excitability   is relatively insensitive to external stimulation,  and the discharge frequency is in a certain range. Meanwhile, from the point of view of bifurcation, type I excitability results from a saddle-node (tangent) bifurcation of equilibrium
points on an invariant circle in ML model, while type II excitability corresponds to a subcritical Hopf bifurcation.
 The ML model is widely used in the theoretical research of the excitatory   nerve discharge \cite{Izhikevich2000spiking,Galan2005,Bogaard2009,Gutkin2005phase,izhikevich2007dynamical,Zhao2016}. Moreover, it can   be used in cardiac cell modeling \cite{zhang2009,yuan2011}.
As we known,  $Na^{+}$ ions and $K^+$ ions cross the  ion channels on the cell membrane back and forth, forming a transmembrane current and leading to the generation of an action potential or a spike - an abrupt and transient change of membrane voltage. An excitable system may be  sensitive to noise \cite{lindner2004effects}.
So an excitable membrane can generate action potential when stimulated by a strong enough input or disturbed by noise.

Recent works on the stochastic ML model  are mostly concerned with  the model under Gaussian noise \cite{Ranjit2017,Montejo2005,li2012spatial,li2014,Jia2015a,Jia2015b}. In order to understand the information coding in the nervous systems, the influence of additive stochastic perturbation on bifurcation scenarios and the stationary distribution of stochastic ML system were studied in \cite{Takashi2004} based on random dynamical systems theory. The Gaussian noise induced multiple spatial coherence resonance and spatial patterns  in  excitable systems  were revealed in \cite{li2012spatial, li2014}. The methods based on asymptotic approximations of the stationary density function and most probable path were developed to understand the role of channel noise in spontaneous excitability \cite{jay2014}. However, Gaussian   noise can not describe some fluctuations with bursts or intermittence or with heavy-tailed distributions, which are   characteristics of $\alpha$-stable L\'evy motions. Indeed,  many complex phenomena involve fluctuations of the L\'evy type, such as asset prices \cite{RM1995economic},   turbulent motions of rotating annular
fluid flows \cite{weeks1995}, a class of biological evolution \cite{GV1996levy}, and random search \cite{Viswanathan2000}. Moreover, recent empirical research has shown that the probability distribution of anomalous (high amplitude) neural oscillations has heavier tail than the standard normal distribution \cite{James2015}.
The neuron systems with L\'evy noise have attracted some  recent attention  \cite{sun2014color,yong2015probability,CR2017FHN,Vinaya2018}.
In fact, L\'evy noise appears to be  more reasonable than Gaussian noise, due to jumps by excitatory and inhibitory impulses caused by external disturbances in biological systems.

In this paper, we will consider the escape problem of the ML model with type II excitability under L\'evy  fluctuations.
In this case the undisturbed system has unique equilibrium state. More concretely, we will study whether the system trajectory starting from the stable equilibrium point in the ML system reaches other region  through a boundary    under the influence of  $\alpha$-stable L\'evy noise. Two different deterministic quantities : the first escape probability (FEP) and the mean first exit time (MFET), are applied to analyse the problem.
In order to quantify  the escape behaviors, we should choose a proper escape region that contains the equilibrium point and a corresponding target region.  The FEP is the likelihood for a system trajectory escaping to the target region, while MFET  is the expected time for a system trajectory exiting the escape region. It turns out that both deterministic quantities are described by nonlocal partial differential equations (details in Sections 3 and 4). Then we will numerically  calculate FEP of the solution stating from the escape region to the target region, and MFET of the solution stating from the escape region to the various outside regions.

The organization of the paper is as follows. In the next section we introduce the undisturbed  ML model and the stochastic
ML model driven by  $\alpha$-stable L\'evy noise.  We also briefly  review  the   $\alpha$-stable L\'evy motion and two deterministic quantities: FEP and MFET, together with   appropriate regions for computing these quantities. In Sections 3 and 4, we  report numerical experiments on the effects of L\'evy motion as well as Brownian motion in ML system, quantified by FEP and by MFET, respectively.
Finally, we end the paper with a summary in Section 5.

\section{The Morris-Lecar model}

We now recall the undisturbed  Morris-Lecar (ML)  model and its disturbed version.

\subsection{The undisturbed  Morris-Lecar model}

The deterministic Morris-Lecar (ML) model has been derived to describe giant barnacle (Balanus Nubilus) muscle fibres \cite{Takashi2004,Tsumoto2006},  and it is represented by the following  two-dimensional system:
\begin{equation}\label{DE}
\begin{split}
C \frac{d{v_{t}}}{dt}&=-g_{Ca}m_{\infty}(v_{t})(v_{t}-V_{Ca})-g_{K}w_{t}(v_{t}-V_{K})-g_{L}(v_{t}-V_{L})+I,\\
\frac{d{w_{t}}}{dt} &=\phi\frac{w_{\infty}(v_{t})-w_{t}}{\tau_{w}(v_{t})},
\end{split}
\end{equation}
where
\begin{align*}
m_{\infty}(v)&= 0.5[1+\tanh(\frac{v-V_{1}}{V_{2}})],\\
w_{\infty}(v)&= 0.5[1+\tanh(\frac{v-V_{3}}{V_{4}})],\\
\tau_{w}(v)&= [\cosh(\frac{v-V_{3}}{2V_{4}})]^{-1}.
\end{align*}
The variables $v_{t}$ and $w_t$ represent the membrane potential and the activation variable for the $K^{+}$ current, respectively. The parameter $C$ stands for the membrane capacitance. The first three terms in the right-hand side of the first equation in the system \eqref{DE} respectively represent the voltage-gated $Ca^{2+}$ current, the voltage-gated delayed rectifier $K^{+}$ current and the leak current. The parameters $g_{Ca}$, $g_K$ and $g_L$ are the maximal conductance of the calcium current, potassium current and leak current, respectively. The parameters $V_{Ca}$, $V_K$ and $V_L$ are the reversal potentials of the calcium current, potassium current and leak current, respectively. Input current is represented by $I$. The constant $\phi$ indicates the change between fast and slow scales of the model. Finally $V_1,~V_2, ~V_3, ~V_4$ are tuning parameters for steady state and time constant.

\begin{figure}[!ht]
\centering 
\includegraphics[height=8cm ,width=10cm]{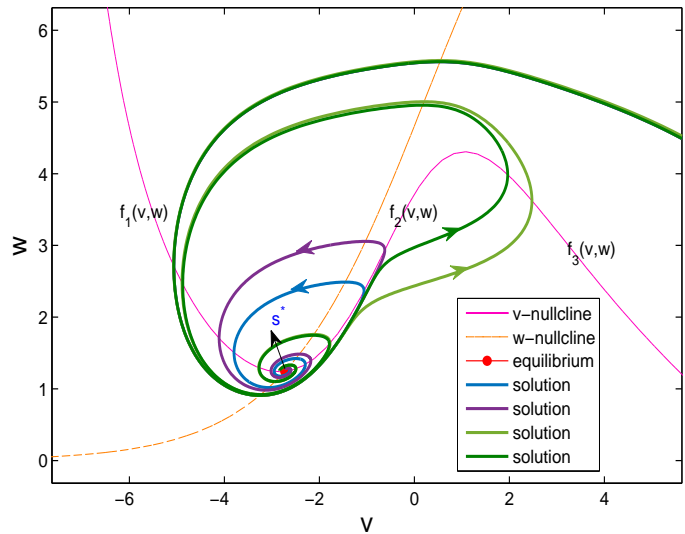}
\caption{The phase portrait of the membrane potential $v_t$ and activation variable $w_t$ in the deterministic Morris-Lecar system \eqref{DE} with type \uppercase\expandafter{\romannumeral2} excitability parameters and $I = 88~(\mu A/cm^2)$.} 
\label{fig:1}
\end{figure}

The parameter values for the type \uppercase\expandafter{\romannumeral2} excitability of ML model are \cite{Rinzel1989,Gutkin1998}: $C = 20~\mu F/cm^2$, $V_{Ca} = 120~mV$, $V_K =
-84~mV$, $V_L = -60~mV$, $g_{Ca} = 4.4~\mu S/cm^2$, $g_K = 8~\mu S/cm^2$, $g_L = 2~\mu S/cm^2$, $V_1 = -1.2~mV$,
$V_2 = 18~mV$, $V_3 = 2~mV$, $V_4 = 30~mV$ and $\phi = 0.04$.  In this case the system possesses a unique
equilibrium state for all values of $I$. This equilibrium is stable for $I<I_H \simeq 93.86~\mu A/cm^2$, and unstable beyond $I_H$
 \cite{Takashi2004,LH2012}. In this study, we choose $I=88~\mu A/cm^2$, so the equilibrium  state is stable.
Furthermore, the resting cells at $I = 88 $ have excitability.
 Figure \ref{fig:1} shows the phase portrait of system \eqref{DE}, in which the $v$-nullcline is divided into three branches, the left branch $f_1$, the middle branch $f_2$ and the right branch $f_3$. Moreover, the resting potential corresponding to the equilibrium $s^*$ is located at the left branch $f_1$.

\subsection{The stochastic Morris-Lecar model}
We consider the Morris-Lecar model driven by symmetric $\alpha$-stable L\'evy motion. This stochastic model is described by the following stochastic differential equations:
\begin{equation}\label{SV}
\begin{split}
&d{v_{t}} =\frac{1}{C}[-g_{Ca}m_{\infty}(v_{t})(v_{t}-V_{Ca})-g_{K}w_{t}(v_{t}-V_{K})-g_{L}(v_{t}-V_{L})+I]dt+\sigma_1 dL_{t}^{1},\\
&d{w_{t}} =\phi\frac{w_{\infty}(x_{t})-w_{t}}{\tau_{w}(v_{t})}dt+\sigma_2 dL_{t}^{2},
\end{split}
\end{equation}
where $L_{t}^{1}$ and $L_{t}^{2}$ are independent scalar symmetric $\alpha$-stable L\'evy motions which have the same jump measure $\nu_{\alpha}$. The symbols $\sigma_1,\sigma_2$ represent
the noise intensities of the $L_{t}^{1}$ and $L_{t}^{2}$,  respectively.
 As a special class of non-Gaussian process with jumps \cite{sato1999levy,bertoin1998levy}, the $\alpha$-stable L\'evy motion is defined by stable L\'evy random variables.
 The distribution for a stable random variable is denoted as $S_{\alpha}(\delta, \beta, \gamma)$.
 Here $\alpha \in (0,2)$ is called the L\'evy motion index (non-Gaussianity index), $\delta$ is the scale parameter, $\beta$ is the skewness parameter, and $\gamma$ is the shift parameter. Let us recall the definition of a symmetric $\alpha$-stable L\'evy motion.

A symmetric $\alpha$-stable L\'evy motion $L_t^\alpha$, with $0 < \alpha < 2$, is a stochastic process with the following properties \cite{applebaum2009levy,duan2015introduction}:

(i) $L_0^\alpha$ = 0, almost surely  (a.s);

(ii) $L_t^\alpha$ has independent increments;

(iii) $L_t^\alpha-L_s^\alpha \sim S_{\alpha}((t-s)^{\frac{1}{\alpha}}, 0, 0)$;

(iv) $L_t^\alpha$ has stochastically continuous sample paths: for every $s$, $L_t^\alpha \to L_s^\alpha$ in probability, as $t \to s$.

The well-known Brownian motion $B_t$ corresponds to  $\alpha$ being $2$.
Moreover, a symmetric $\alpha$-stable L\'evy motion can be represented as the triplet (0,~0,~$\nu_{\alpha}$), where the jump measure $\nu_{\alpha}$ is defined as \cite{sato1999levy,samorodnitsky1996stable}
\begin{align}
  \nu_{\alpha}=\frac{C_{\alpha}dy}{| y | ^{1+\alpha}},
\end{align}
with
\begin{align}
 C_{\alpha}=\frac{\alpha}{2^{1-\alpha}\pi}\frac{\Gamma(1+\frac{\alpha}{2})}{\Gamma(1-\frac{\alpha}{2})}.
\end{align}

For $0<\alpha<2$, the following tail estimate of stable L\'evy random variable $L$ holds \cite{samorodnitsky1996stable}
\begin{align}
  \lim_{y\rightarrow\infty}y^{\alpha}\mathbb{P}(L>y)=C_{\alpha}\frac{1+\beta}{2}\sigma^{\alpha}.
\end{align}
This estimate indicates that the stable L\'evy random variable $L$ has a ``heavy tail", which decays polynomially, unlike the tail estimate of Gaussian random variable, which decays exponentially.

In this paper, we make scale transformation of variables for the convenience of calculation. After this transformation, the stable equilibrium is $s^* = (-2.7277, 1.2436)$ in the deterministic case (see equation \eqref{DE}).
The noise term $L_{t}^{1}$ represents current fluctuations in the environment and the noise term $L_{t}^{2}$ is understood as ion channel noise due to the opening and closing of ion channels. Indeed, there exists small fluctuations, always combined with  unpredictable jumps of random environments in the biology growth process.
 L\'evy motion $L_t^{\alpha}$ is suitable for simulating this type of noise. Moreover, $\alpha$-stable L\'evy motion has larger jumps with lower jump frequencies when $\alpha$ closes to 0, while for $1 < \alpha < 2$, it has smaller jumps with higher jump probabilities. Especially, in the following numerical simulations, when $\alpha$ being 2, we respectively change the L\'evy motion $L_t^{1}$, $L_t^{2}$ to independent Brownian motion $B_t^1$, $B_t^2$.

\begin{figure}[!ht]
\begin{minipage}{0.49 \textwidth}
\centerline{\includegraphics[width=9cm,height=8cm]{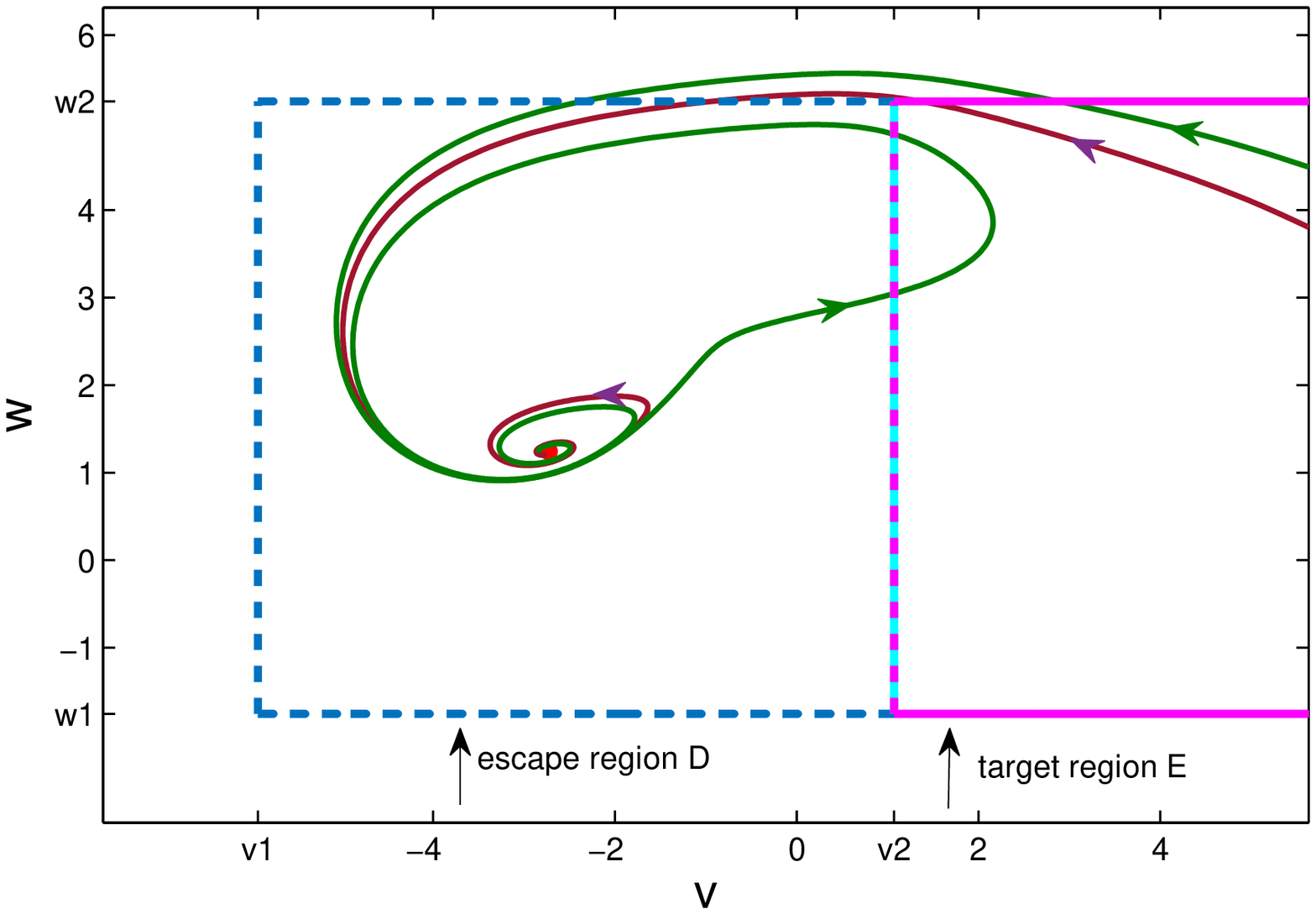}}
\centerline{(a)}
\end{minipage}
\hfill
\begin{minipage}{0.49 \textwidth}
\centerline{\includegraphics[width=9cm,height=8cm]{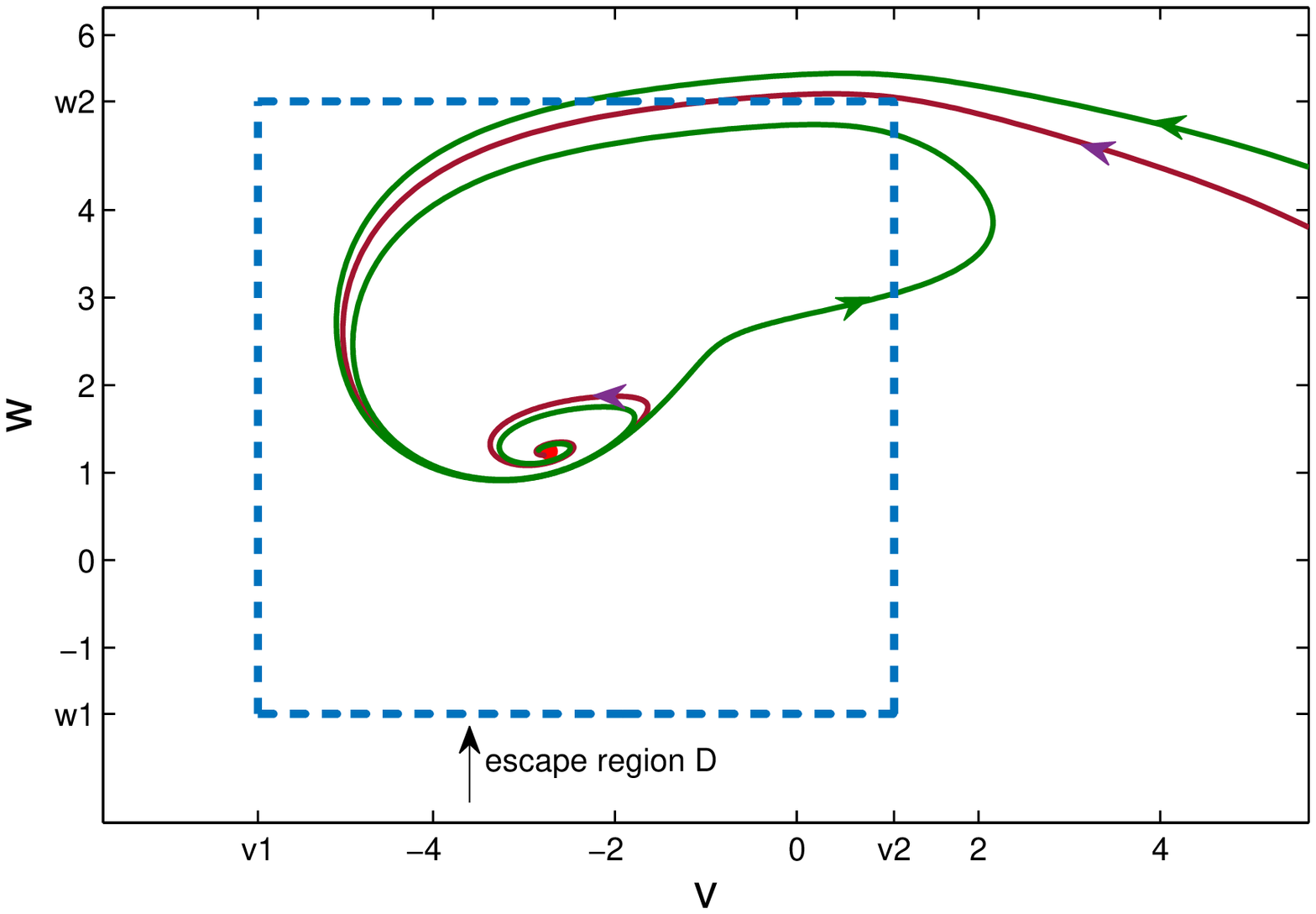}}
\centerline{(b)}
\end{minipage}
\vfill
\caption{(a) The escape region we choose for   FEP is $D:(-5.9277,1.0723)\times(-1.7564,5.2436)$ and target region  is $E: [1.0723,\infty)\times[-1.7564,5.2436]$. (b) The escape region $D$ for   MFET  is the same as in (a), while the  target region is $D^{c}$.}
\label{fig:2}
\end{figure}

 If ML system \eqref{DE} is not subject to any disturbance, the stable equilibrium corresponds to the resting state of the ML system. When the system is under a perturbation, the solution orbit (or path) starting from the resting state may  respond  as a small oscillation near the resting state or produce a spike. This means that the ML system may have a state transition under the interference of noise. We further note that the $v$-nullcline is ``cubic'', the middle branch $f_2$ in some sense separates the firing of an action potential from the subthreshold return to equilibrium \cite{Ermentrout2010}. If an orbit crosses the separatrix, it will be attracted by the right branch $f_3$ under the interference of non-Gaussian noise. Then we explain whether the orbit from the equilibrium state can be attracted by $f_3$ as an escape problem.
 Two deterministic quantities: FEP and MFET, are applied to this problem. In order to calculate FEP of the equilibrium $s^*$, the escape region $D$ (containing $s^*$) and the target region $E$ should be chosen. The two
regions are $D:(-5.9277,1.0723)\times(-1.7564,5.2436)$ and $E: [1.0723,\infty)\times[-1.7564,5.2436]$
as shown in Figure \ref{fig:2}. The FEP represents the probability of the solution orbit starting at a point in $D$ first escapes to region $E$. The reason we choose $\{(v,w)|v=1.0723, w \in R \}$ as a boundary line to calculate FEP is that $v > 1.07234$ can be seen as the high  electrical potential of a nerve cell. And if the solution orbit crosses this line $\{(v, w)|v=1.0723, w \in R \}$, it must be attracted by $f_3$. Meanwhile, the stability of the equilibrium under the stimulation by  L\'evy Motion and Brownain Motion will also be considered. We select region $D$ to compute MFET, which implies the mean first exit time of the solution orbit starting at a point in region $D$ escapes to region $D^{c}$.

\section{First escape probability}

The general form of the two-dimensional stochastic differential system \eqref{SV} is as follows:
\begin{equation}\label{general}
\begin{split}
&d{v_t} = f_1(v_t,w_t)dt  + \sigma_1 dL_{t}^{1},\\
&d{w_t} =f_2(v_t,w_t)dt  +\sigma_2 dL_{t}^{2}.
\end{split}
\end{equation}
The infinitesimal generator $A$ of the system \eqref{general} is
\begin{align}\label{GNR}
Ap(v,w)=&f_{1}(v,w)p_{v}(v,w)+f_{2}(v,w)p_{w}(v,w)  \nonumber\\
& +\sigma_{1}^{\alpha}\int_{\mathbb{R}\setminus{\{0\}}}[p(v+v',w)-p(v,w)]\nu_{\alpha}(dv')
\nonumber\\
&+\sigma_{2}^{\alpha}\int_{\mathbb{R}\setminus{\{0}\}}[p(v,w+w')-p(v,w)]\nu_{\alpha}(dw').
\end{align}
When  $L_{t}^{1}, L_{t}^{2}$  are replaced by independent Brownian motions $B_t^1, B_t^2$, the generator becomes
\begin{equation}\label{GNRbm}
Ap(v,w)=f_{1}(v,w)p_{v}(v,w)+f_{2}(v,w))p_{w}(v,w)+\frac{\sigma_{1}^{2}}{2}p_{vv}(v,w)
+\frac{\sigma_{2}^{2}}{2}p_{ww}(v,w).
\end{equation}

\begin{figure}
\begin{minipage}[!h]{0.49 \textwidth}
\centerline{\includegraphics[width=7cm,height=3.8cm]{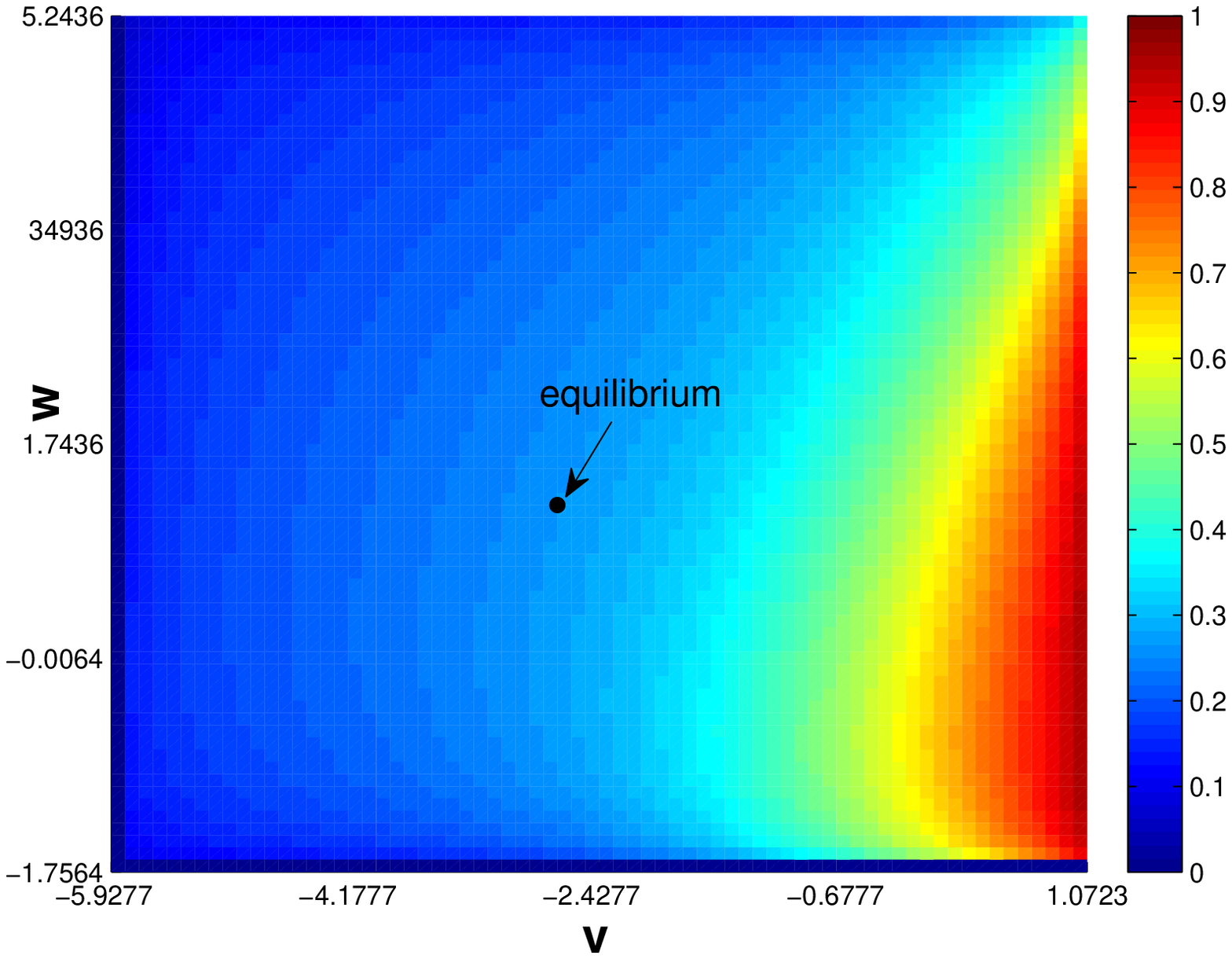}}
\centerline{(a) $\alpha$ = 0.5, $\sigma$ = 0.5.}
\end{minipage}
\hfill
\begin{minipage}[!h]{0.49 \textwidth}
\centerline{\includegraphics[width=7cm,height=3.8cm]{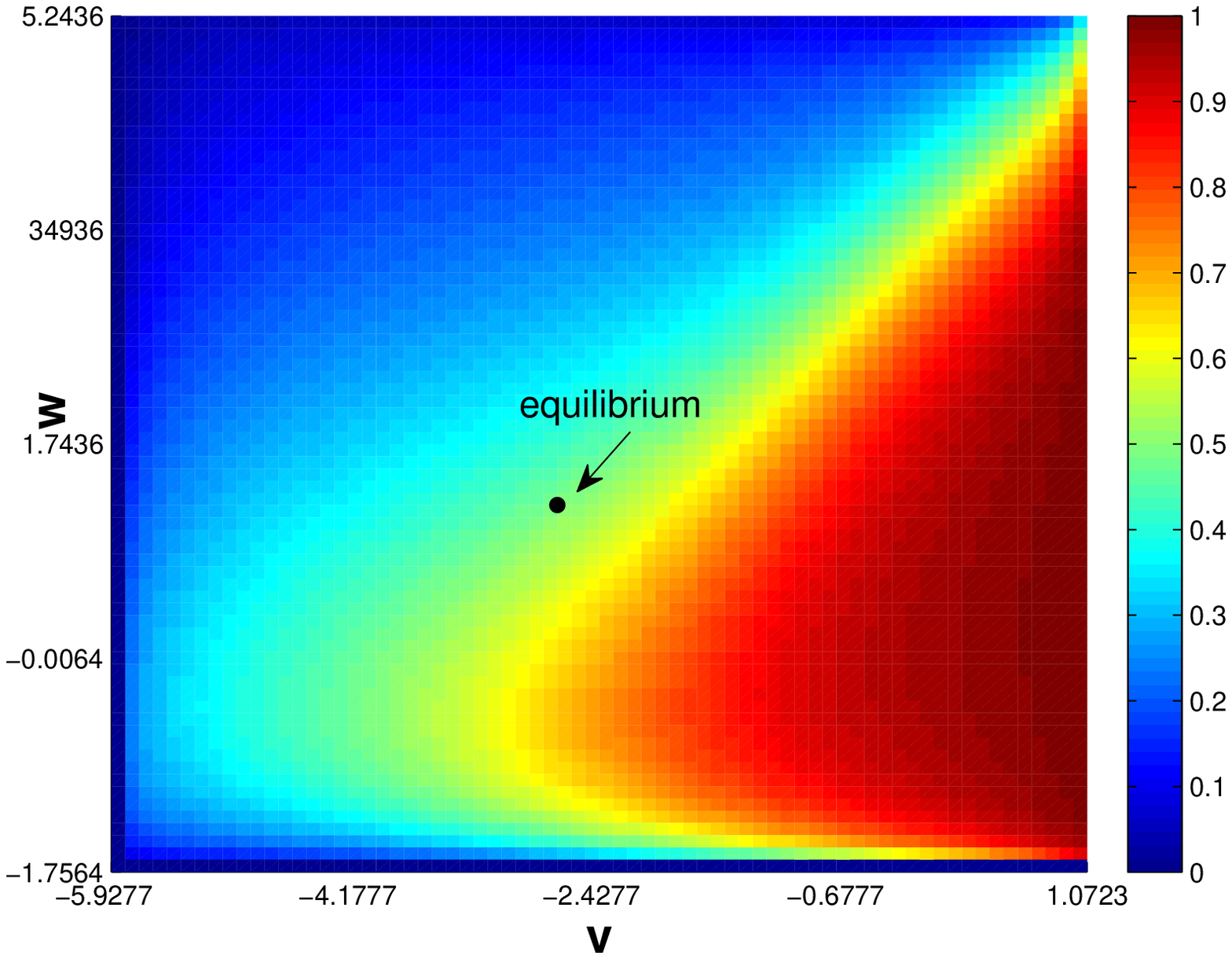}}
\centerline{(e) $\alpha$ = 1.25, $\sigma$ = 0.25.}
\end{minipage}

\begin{minipage}[!h]{0.49 \textwidth}
\centerline{\includegraphics[width=7cm,height=3.8cm]{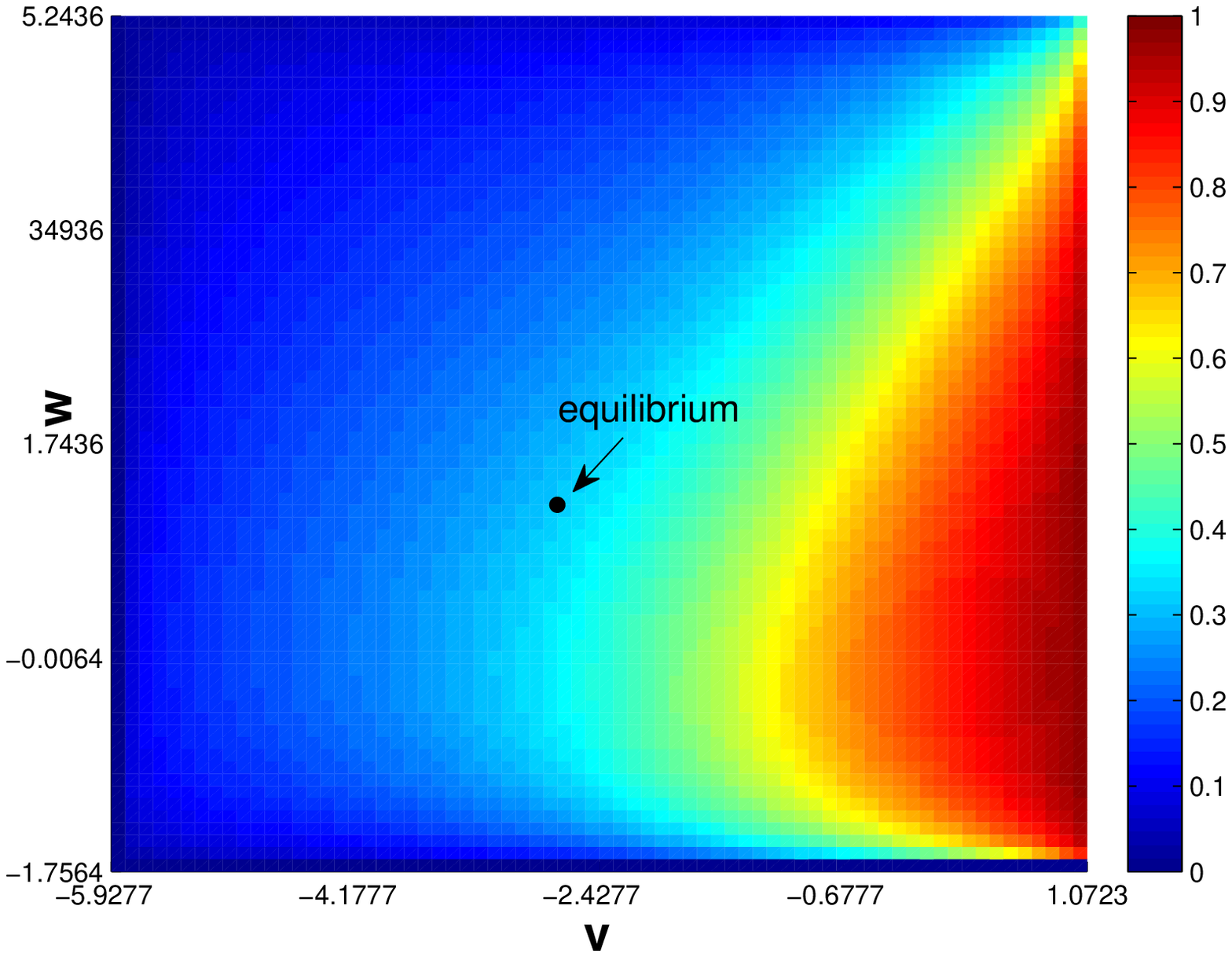}}
\centerline{(b) $\alpha$ = 1, $\sigma$ = 0.5.}
\end{minipage}
\hfill
\begin{minipage}[!h]{0.49 \textwidth}
\centerline{\includegraphics[width=7cm,height=3.8cm]{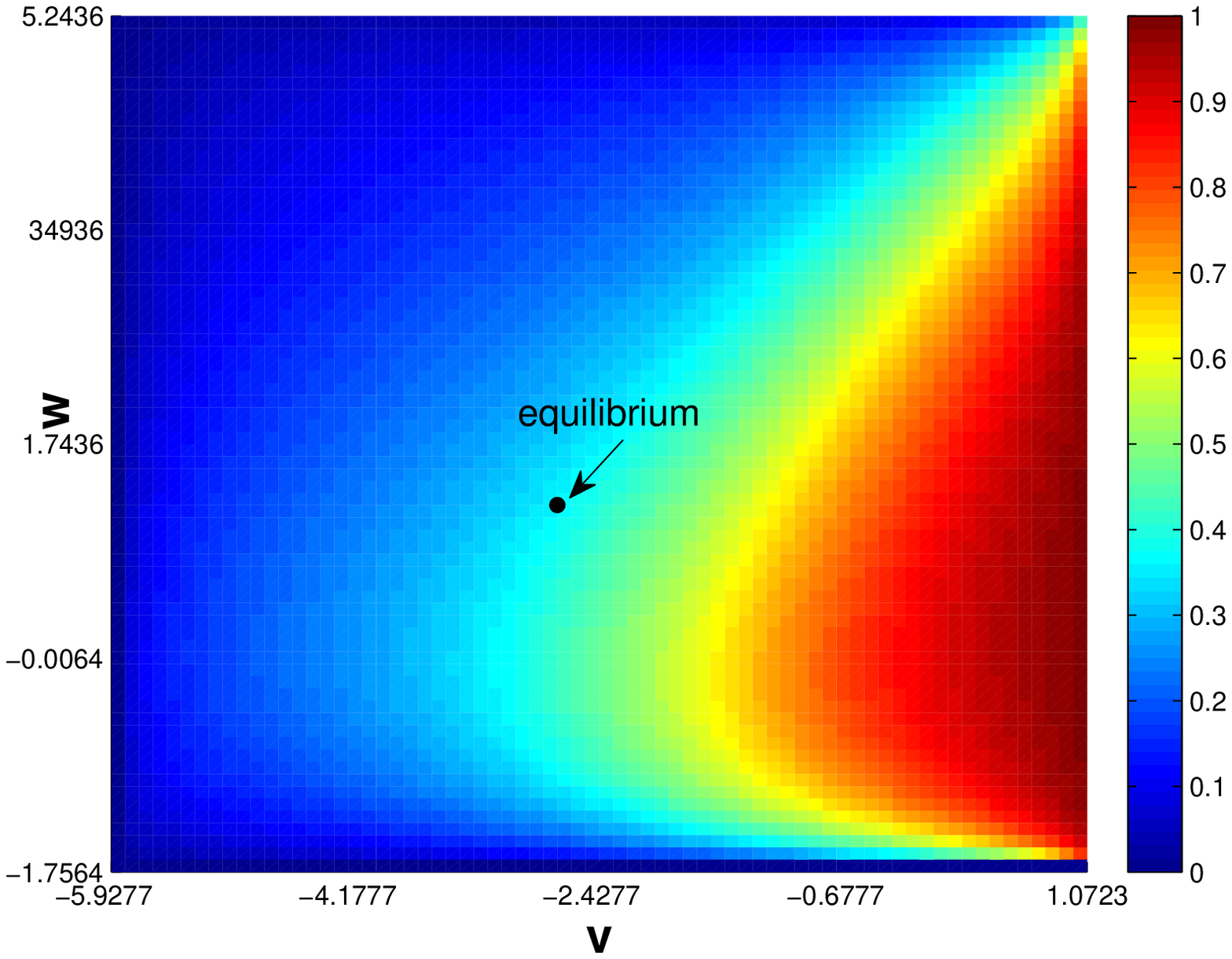}}
\centerline{(f) $\alpha$ = 1.25, $\sigma$ = 0.5.}
\end{minipage}

\begin{minipage}[!h]{0.49 \textwidth}
\centerline{\includegraphics[width=7cm,height=3.8cm]{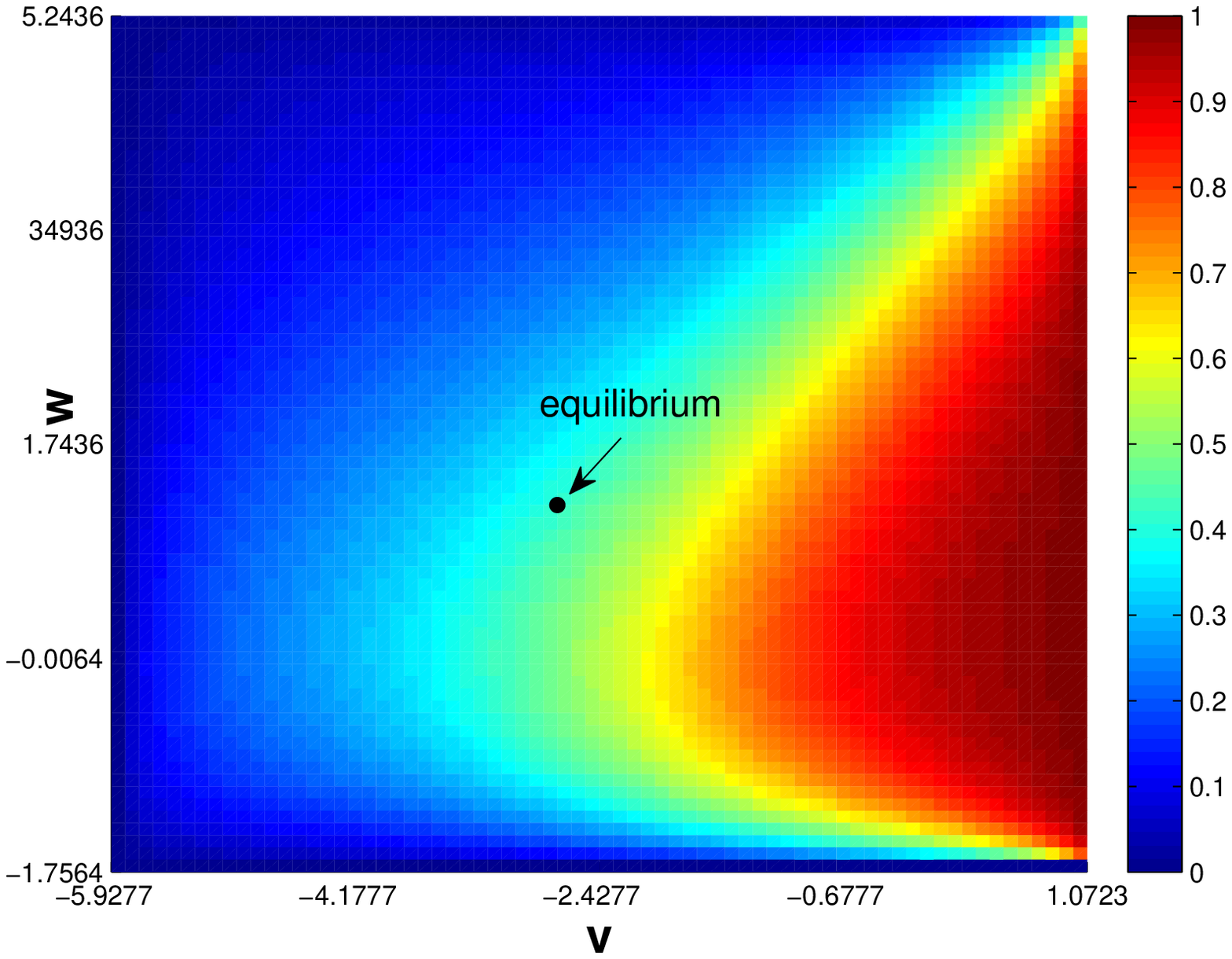}}
\centerline{(c) $\alpha$ = 1.5, $\sigma$ = 0.5.}
\end{minipage}
\hfill
\begin{minipage}[!h]{0.49 \textwidth}
\centerline{\includegraphics[width=7cm,height=3.8cm]{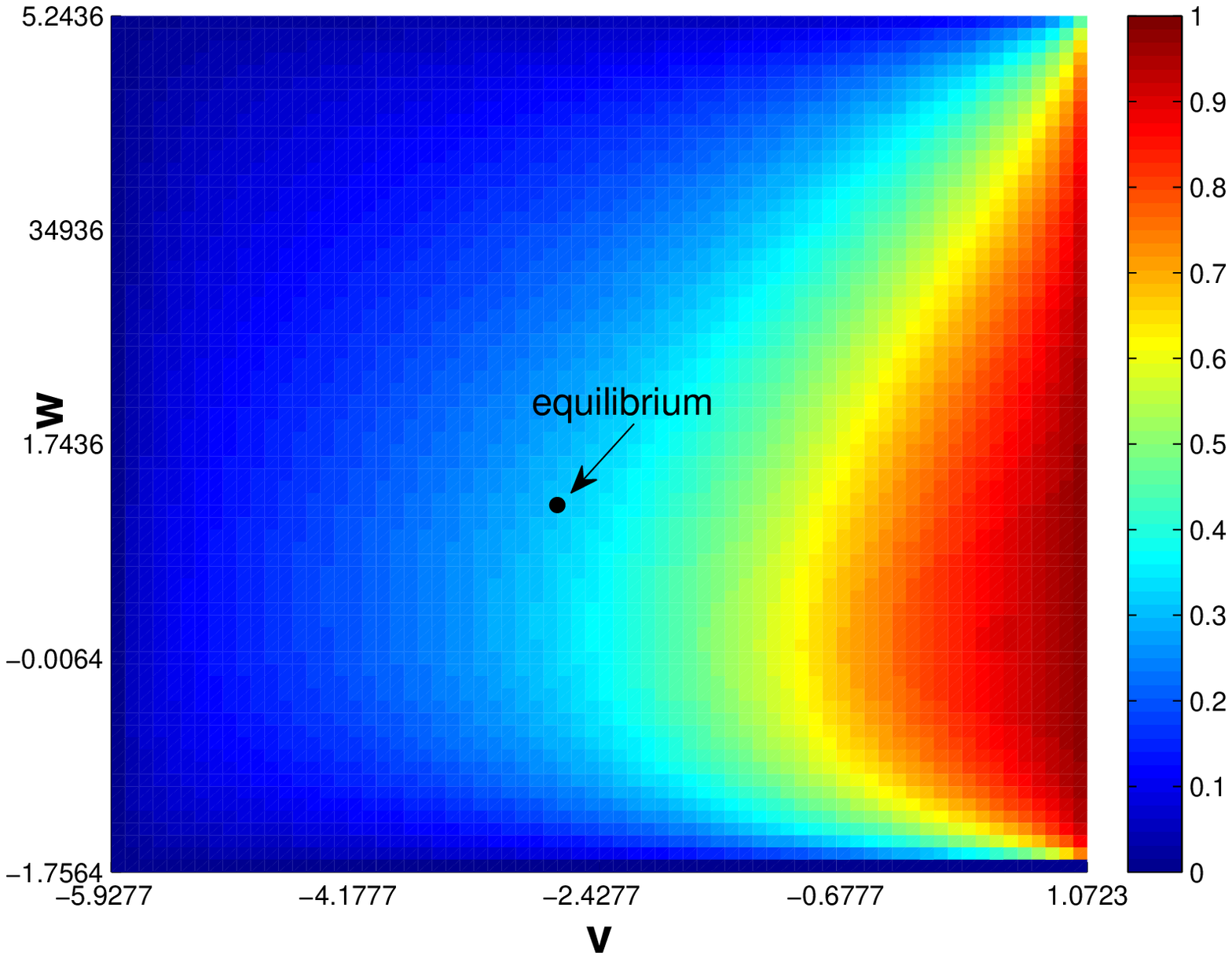}}
\centerline{(g) $\alpha$ = 1.25, $\sigma$ = 0.75.}
\end{minipage}

\begin{minipage}[!h]{0.49 \textwidth}
\centerline{\includegraphics[width=7cm,height=3.8cm]{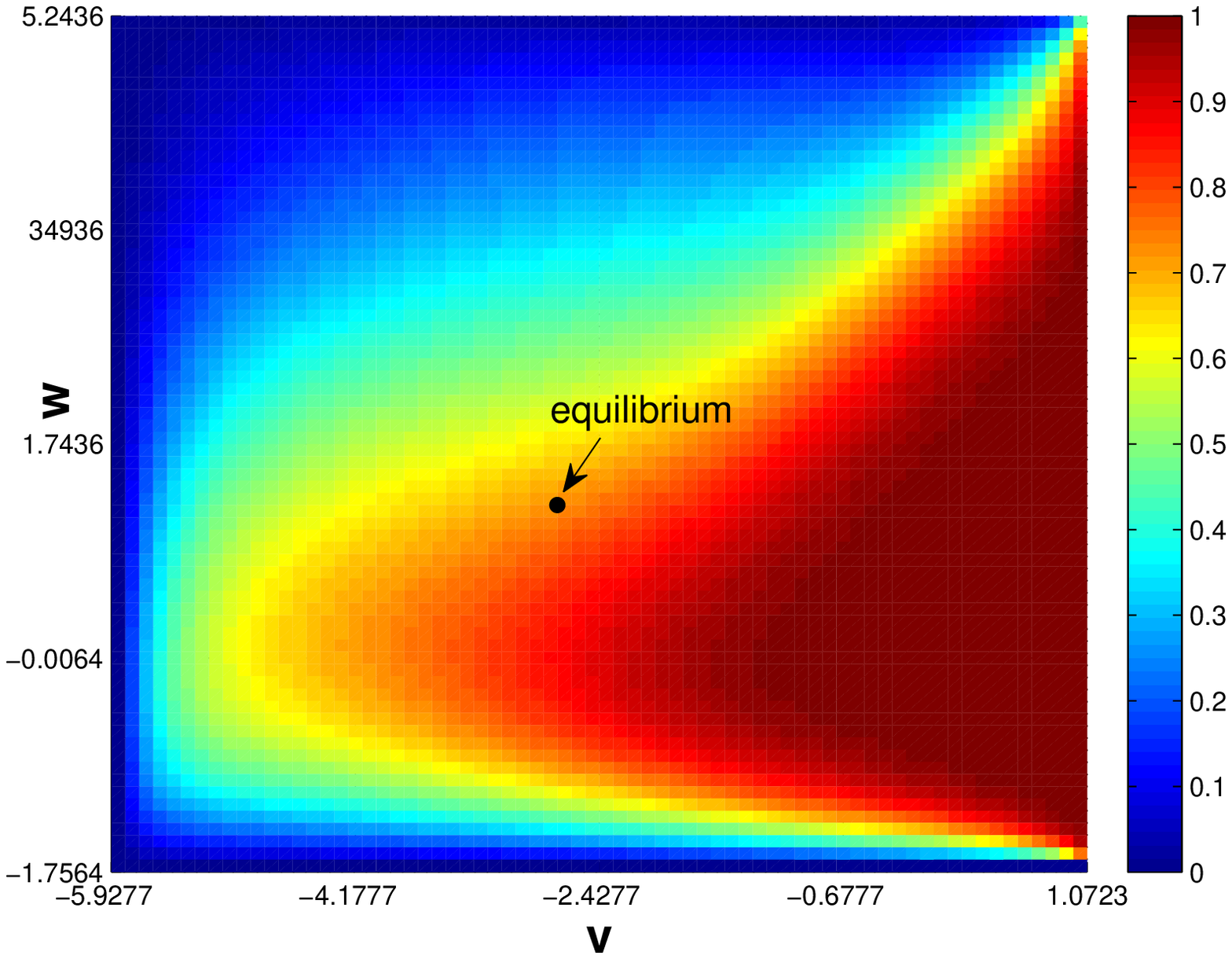}}
\centerline{(d) Brownain case, $\sigma$ = 0.5.}
\end{minipage}
\hfill
\begin{minipage}[!h]{0.49 \textwidth}
\centerline{\includegraphics[width=7cm,height=3.8cm]{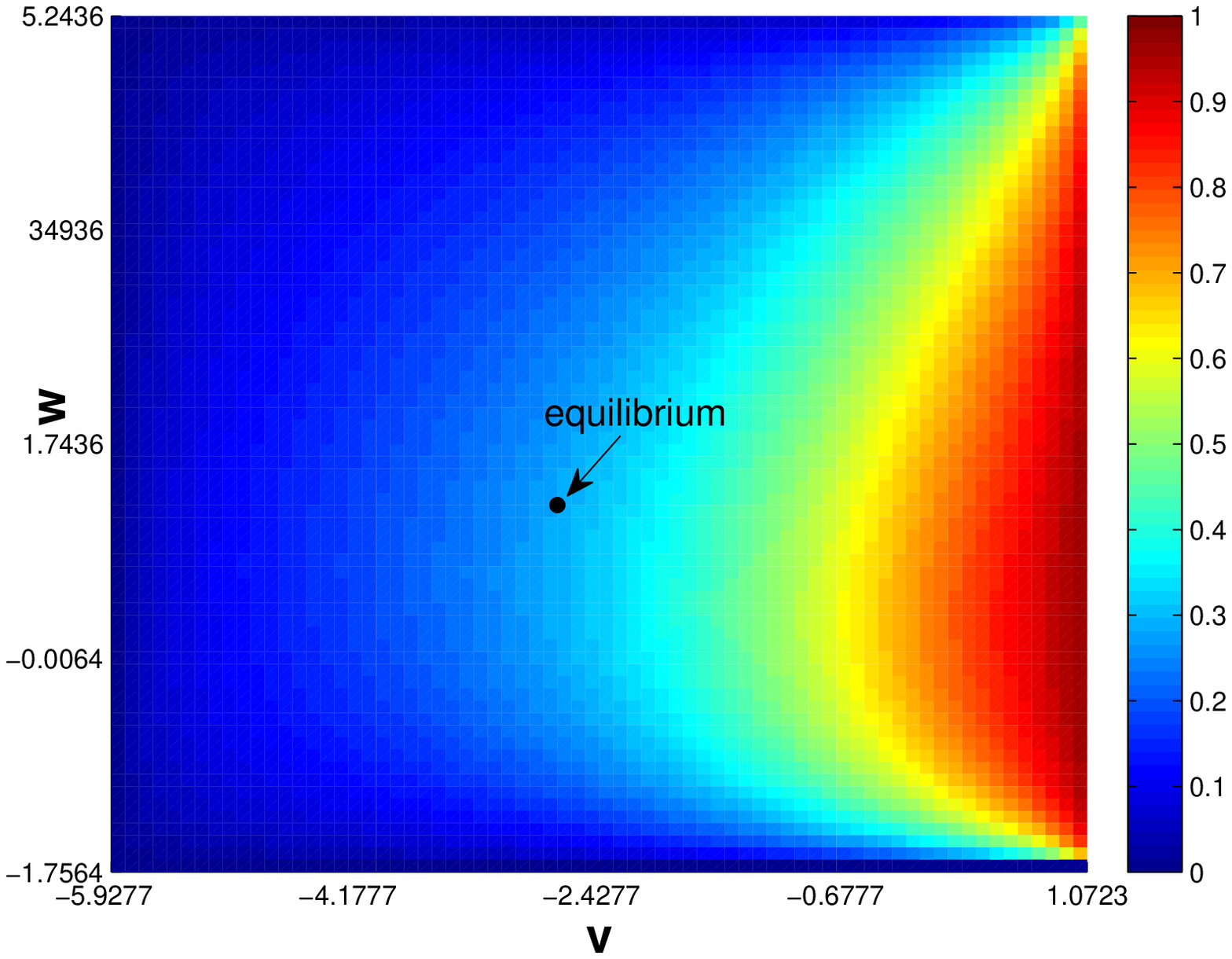}}
\centerline{(h) $\alpha$ = 1.25, $\sigma$ = 1.}
\end{minipage}

\vfill
\caption{FEP $p(v,w)$ from the escape region $D:(-5.9277,1.0723)\times(-1.7564,5.2436)$ to target region
$E: [1.0723,\infty)\times[-1.7564,5.2436]$. The color map depends on L\'evy motion index $\alpha$ and noise intensity $\sigma$ ($\sigma_1= \sigma_2=\sigma$). (a)-(d) Influence of L\'evy motion index $\alpha$ on FEP for different values of $\alpha$ with fixed noise intensity $\sigma = 0.5$. (e)-(h) Influence of noise intensity $\sigma$ on FEP for different values of $\sigma$ with fixed L\'evy motion index $\alpha=1.25$. The color bar in all figures is set to the same scale, red making 1 and blue making 0.}
\label{fig:3}
\end{figure}

\begin{figure}[!ht]
\begin{minipage}{0.49 \textwidth}
\centerline{\includegraphics[width=9cm,height=6cm]{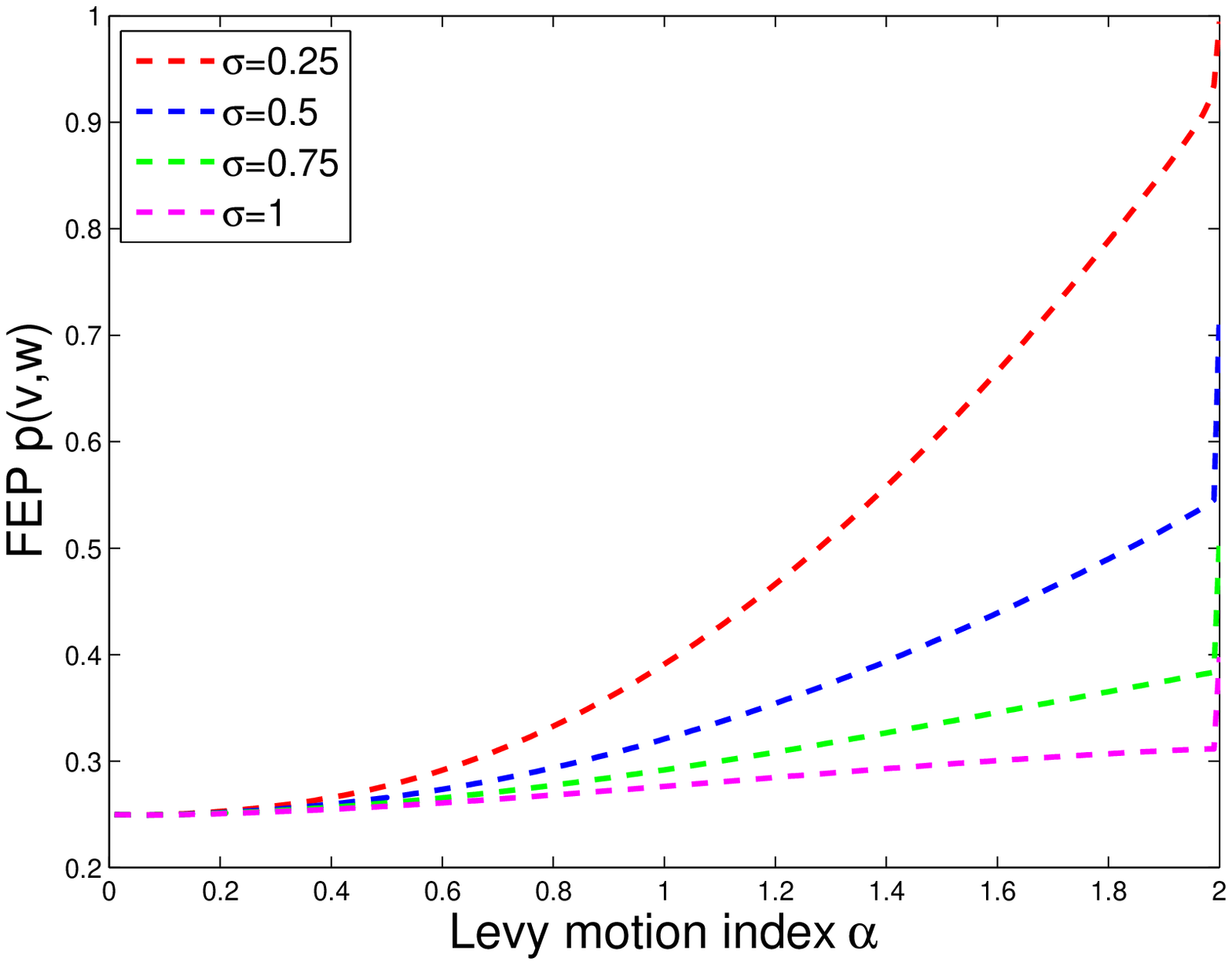}}
\centerline{(a)}
\end{minipage}
\hfill
\begin{minipage}{0.49 \textwidth}
\centerline{\includegraphics[width=9cm,height=6cm]{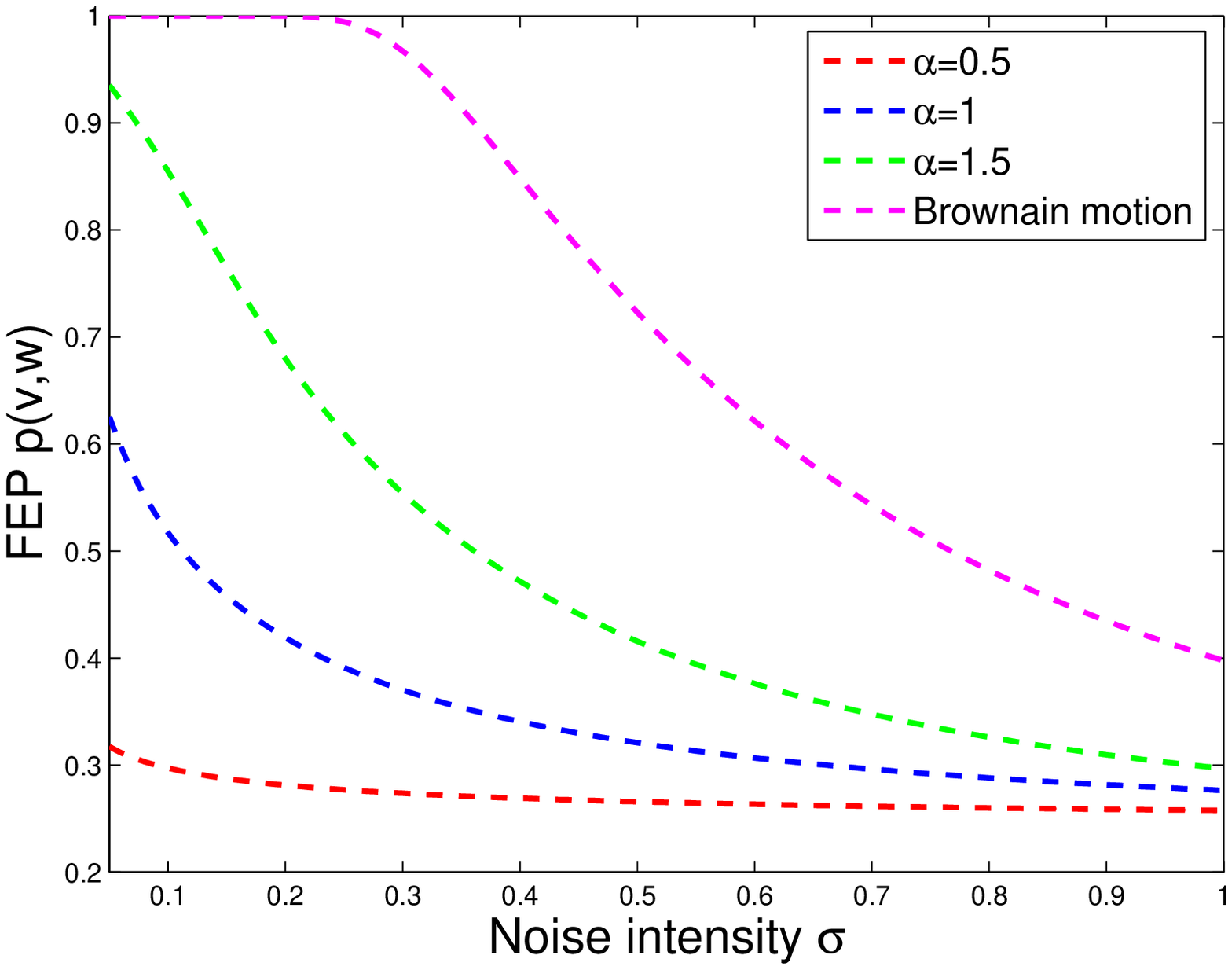}}
\centerline{(b)}
\end{minipage}

\vfill
\caption{FEP $p(v,w)$ at the equilibrium $s^* = (-2.7277, 1.2436)$, $\sigma_1=\sigma_2=\sigma$. (a) Effect of L{\'e}vy motion index $\alpha$ on FEP at $s^*$ for various values of noise intensity $\sigma$ (red: $\sigma$ = 0.25, blue: $\sigma$  = 0.5, green: $\sigma$  = 0.75, pink: $\sigma$  = 1). (b) Effect of noise intensity $\sigma$ on FEP at $s^*$ for various values of L\'evy motion index $\alpha$ (red: $\alpha$ = 0.5, blue: $\alpha$  = 1, green: $\alpha$  = 1.5, pink: $\alpha$ being 2 indicates to case of Brownain motion).}
\label{fig:4}
\end{figure}

The first escape probability (FEP) here is employed to characterize the likelihood of the escape   from the low potential region to the high potential region for the ML system. It is denoted as $p(v,w)$, which is used to characterize the probability of the solution orbit starting at $(v,w)$ in a open region $D$ first escaping to a target region $E$.
The escape probability $p(v,w)$ can be solved by the following  integral-differential equation with a Balayage-Dirichlet exterior boundary condition  \cite{duan2015introduction}:
\begin{align}\label{FEP}
  &Ap(v,w)=0, \quad (v,w)\in D,\notag \\
  &p(v,w)=
  \begin{cases}
    1, \quad (v,w)\in E,  \\
    0, \quad (v,w)\in D^{c}\setminus{E}.
  \end{cases}
\end{align}
Here $ D^{c}$ is the complement set of the bounded region $D$. Results about the existence and uniqueness of solutions to nonlocal systems similar to integral-differential equation \eqref{FEP} may be found in \cite{MG2002,TK2009,DD2019}. The equation \eqref{FEP} can be solved by an effective numerical scheme given in the Appendix.

Figure \ref{fig:3} shows the numerical simulation of the FEP for different L\'evy motion index $\alpha$ and noise intensity $\sigma$ ($\sigma_1=\sigma_2=\sigma$). In Figure \ref{fig:3}(a)-(d), the noise intensity $\sigma=0.5$, the area of high FEP (red region) gradually becomes bigger with the increase of $\alpha$.
 In Figure \ref{fig:3}(e)-(h), when the L\'evy motion index $\alpha=1.25$, the area of high FEP (red region) gradually becomes smaller with the increase of $\sigma$.  From Figure \ref{fig:3}(a)-(h), we conclude that the escape probability varies with initial membranae potential, in the following way: when the noise intensity is fixed ($\sigma=0.5$), the larger L\'evy motion index $\alpha$, the more likely the low potential state of the ML system becomes excited. However, when the  L\'evy motion index $\alpha$ is fixed ($\alpha=1.25$), the smaller the noise intensity, the easier the  ML system gets excited.

Figure \ref{fig:4} depicts the influence of  L\'evy motion index $\alpha$ and  noise intensity $\sigma$ ($\sigma_1=\sigma_2=\sigma$) on FEP of the equilibrium $s^*$.
In Figure \ref{fig:4}(a),
 L\'evy motion index $\alpha$ is from 0.01 to 1.99 and the step size is 0.01, the value of $\alpha$ being 2 corresponds to the Brownain Motion case. It can be seen that for fixed noise intensity $\sigma$, FEP increases as $\alpha$ increases, and grows rapidly at $\alpha$ being 2.
 This means that smaller jumps with higher frequencies are more beneficial to  make the system from  the resting state  to the excited state and generate a pulse. In Figure \ref{fig:4}(b), for various value of fixed $\alpha$ ($\alpha$ = 0.5, 1, 1.5, 2), the smaller noise intensity, the higher FEP,  the easier the system gets into the excited state.
 When $\alpha=0.5$, FEP keeps almost unchanged as the noise intensity increases, which means that noise has less effect on the state transitions.
  When $\alpha=1$ and $\alpha=1.5$, FEP drops rapidly with small noise intensity and as the noise intensity increases, the trend of FEP tends to be steady.
 Finally, when $\alpha$ being 2, with increasing $\sigma$ from 0.05, FEP remains at 1 until $\sigma_a$ (=0.185),
 and then it becomes monotonously decreasing. This means when the system is affected by the Brownian motion, the solution path starting at the equilibrium point can definitely escape region $D$ for very small noise intensity.
   In general, the larger the L\'evy motion index $\alpha$ and the smaller noise intensity $\sigma$, the more likely for the ML system to reach the excited state and generate a pulse. Compared with non-Gaussian L\'evy noise, the Gaussian Brownian noise is easier to make the system transition from the resting state to the excited state.

\begin{figure}[!ht]
\begin{minipage}{0.49 \textwidth}
\centerline{\includegraphics[width=9cm,height=6cm]{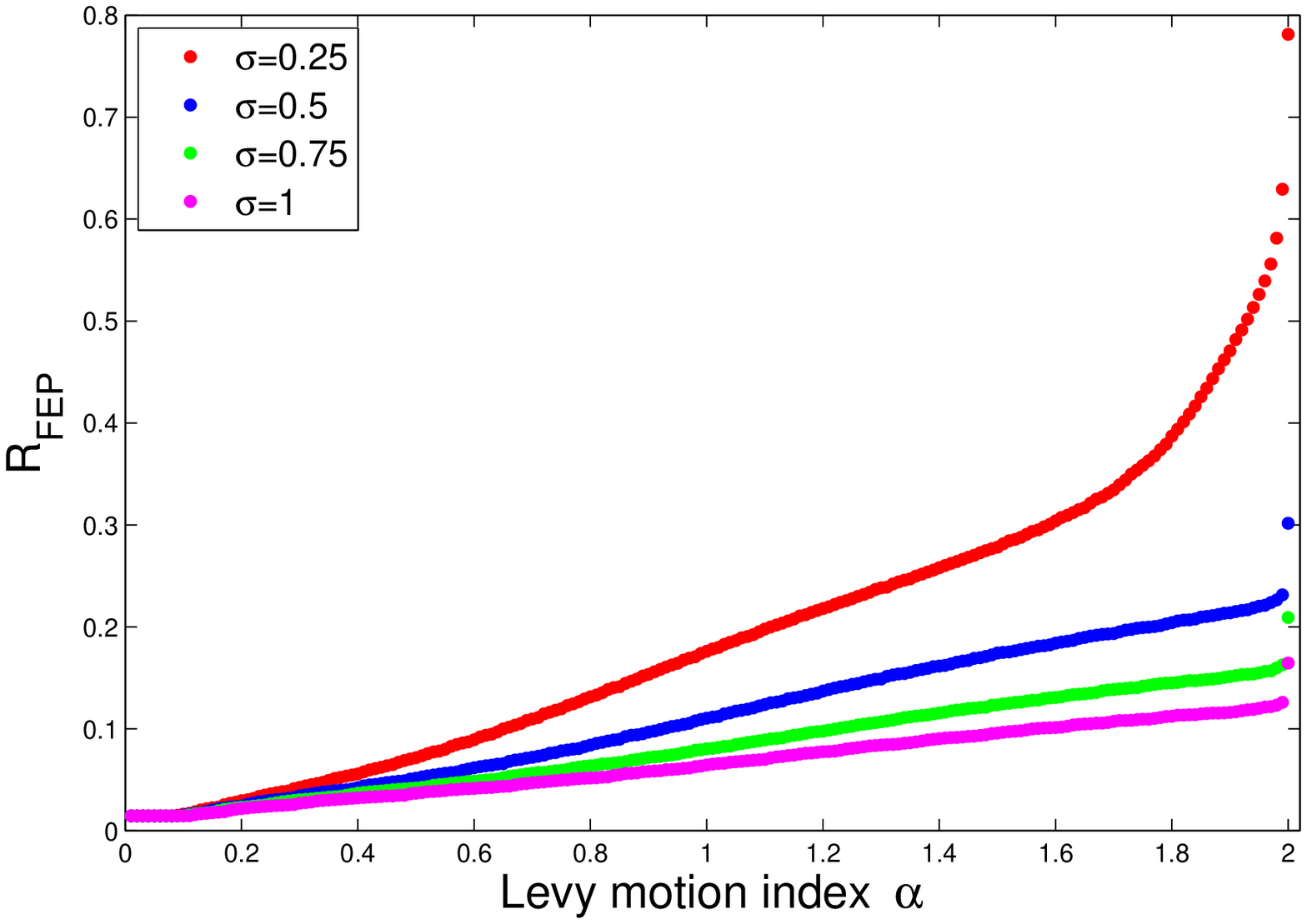}}
\centerline{(a)}
\end{minipage}
\hfill
\begin{minipage}{0.49 \textwidth}
\centerline{\includegraphics[width=9cm,height=6cm]{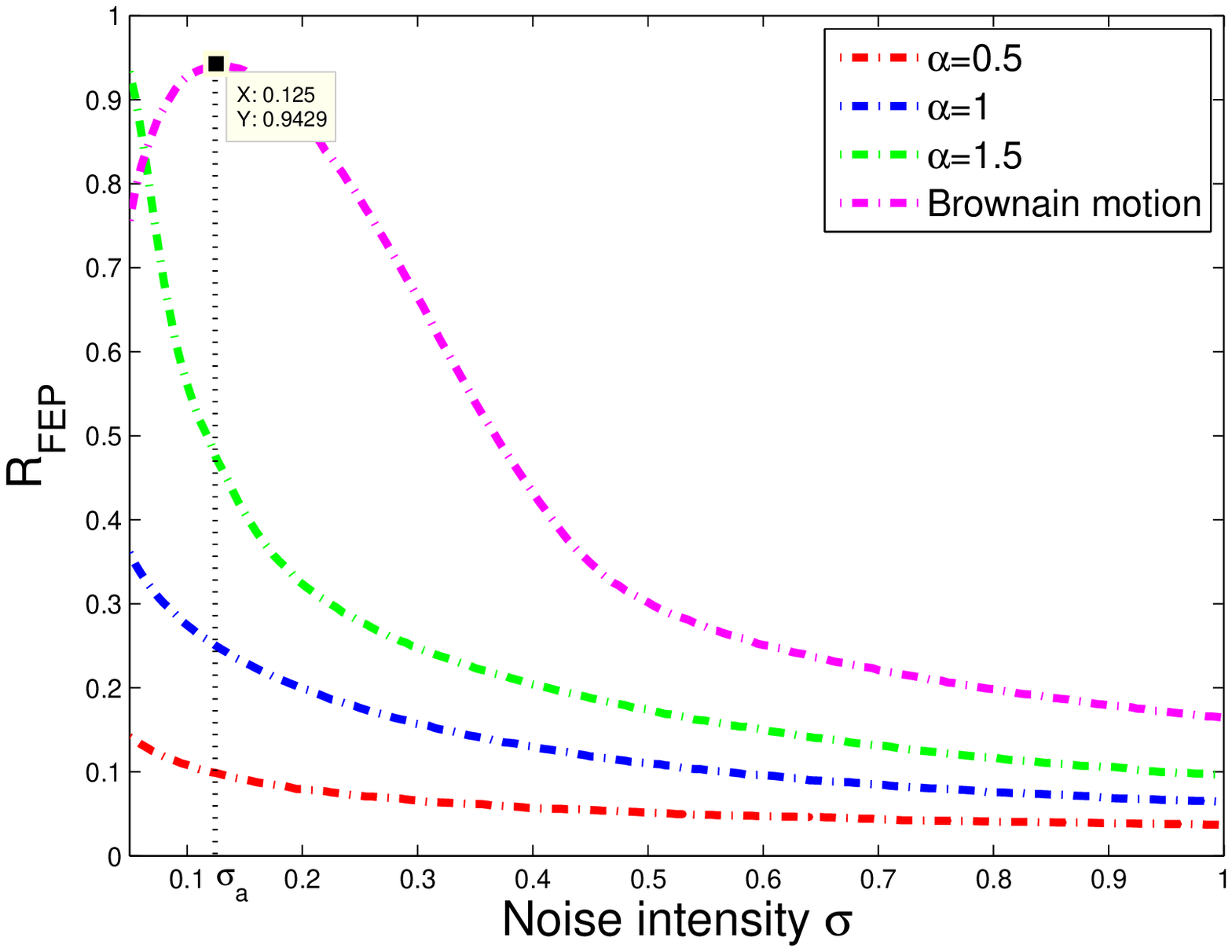}}
\centerline{(b)}
\end{minipage}

\vfill
\caption{Effect of noise intensity and L{\'e}vy motion index on the $R_{FEP}$. The
threshold $p^*$ to compute $R_{FEP}$ is chosen as 0.8. (a) $R_{FEP}$ against L\'evy motion
index $\alpha$ for various noise intensity $\sigma$ (red: $\sigma$ = 0.25, blue: $\sigma$ = 0.5, green: $\sigma$ = 0.75,
pink: $\sigma$ = 1). (b) $R_{FEP}$ against noise intensity $\sigma$ for various L{\'e}vy motion index $\alpha$ (red: $\alpha$ = 0.5, blue: $\alpha$ = 1, green: $\alpha$ = 1.5, pink: $\alpha$ being 2 indicates to case of Brownain motion.}
\label{fig:5}
\end{figure}

\begin{figure}[!ht]
\begin{minipage}{0.49 \textwidth}
\centerline{\includegraphics[width=9cm,height=6cm]{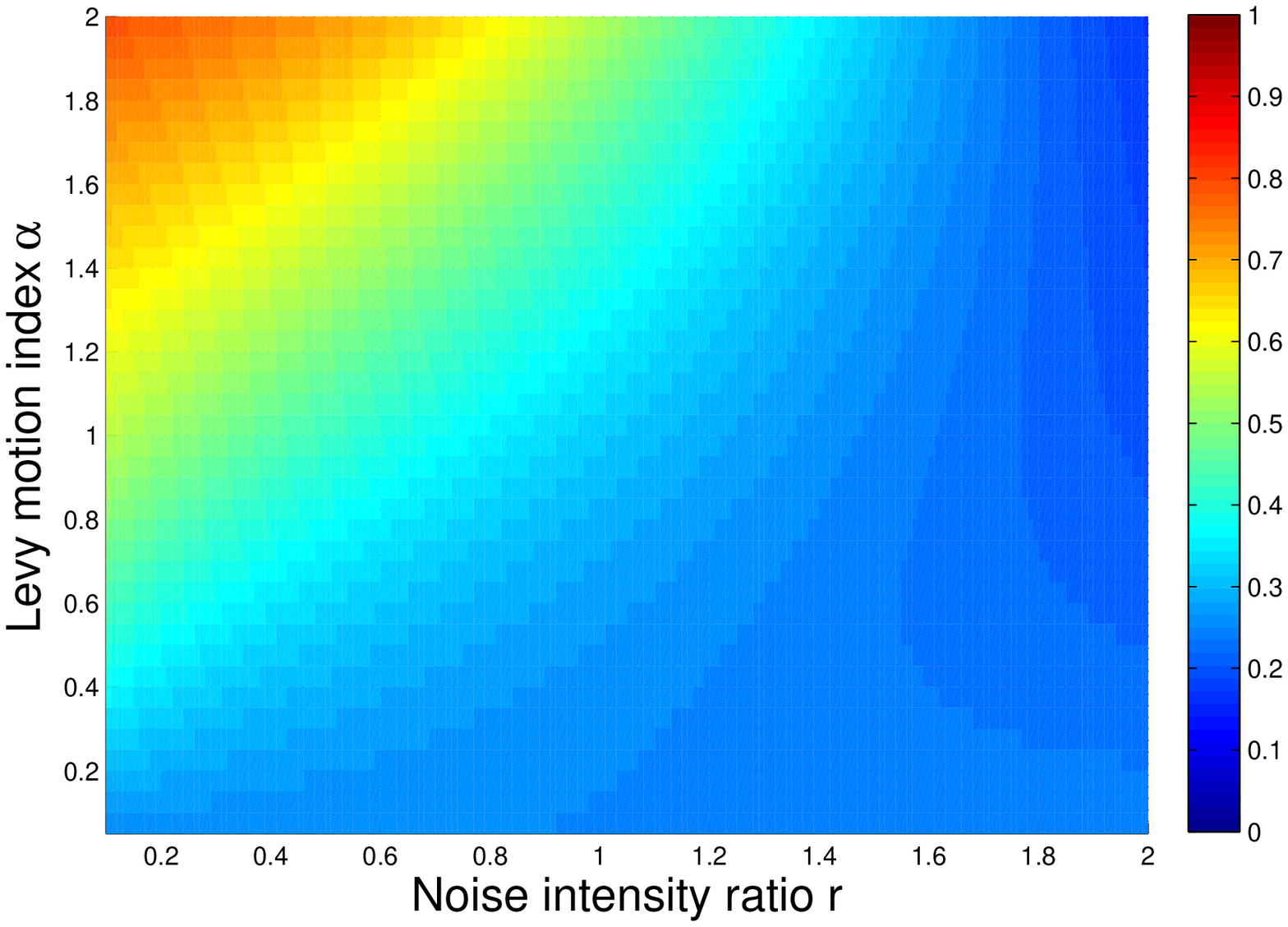}}
\centerline{(a) $r=\sigma_2/\sigma_1$, fixed $\sigma_1=0.5$.}
\end{minipage}
\hfill
\begin{minipage}{0.49 \textwidth}
\centerline{\includegraphics[width=9cm,height=6cm]{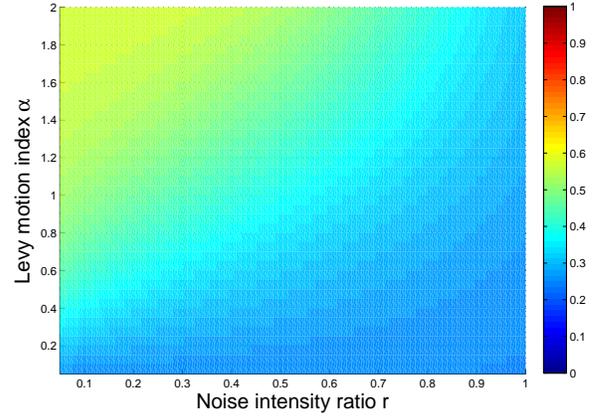}}
\centerline{(b) $r=\sigma_2/\sigma_1$, fixed $\sigma_1=1$.}
\end{minipage}

\begin{minipage}{0.49 \textwidth}
\centerline{\includegraphics[width=9cm,height=6cm]{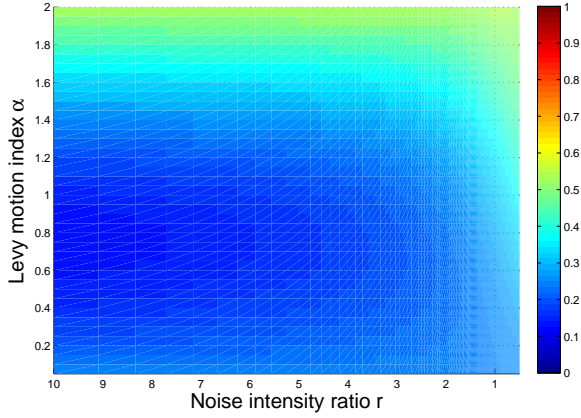}}
\centerline{(c) $r=\sigma_2/\sigma_1$, fixed $\sigma_2=0.5$.}
\end{minipage}
\hfill
\begin{minipage}{0.49 \textwidth}
\centerline{\includegraphics[width=9cm,height=6cm]{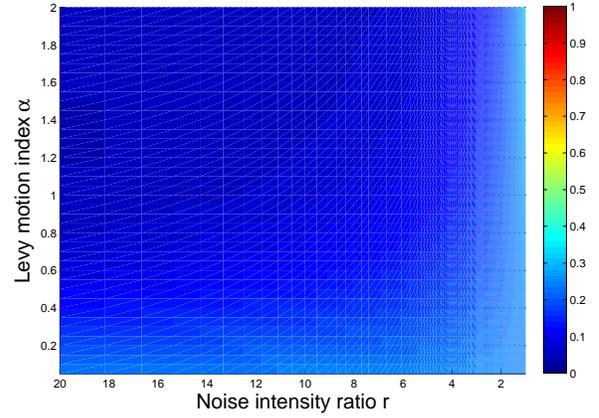}}
\centerline{(d) $r=\sigma_2/\sigma_1$, fixed $\sigma_2=1$.}
\end{minipage}

\vfill
\caption{FEP $p(v,w)$ at the equilibrium $s^*$ = (-2.7277, 1.2436), the color map depends on  noise intensity ratio $r$ and L\'evy motion index $\alpha$. The color bar in all figures is set to the same scale, red making 1 and blue making 0. (a) Fixed $\sigma_1 = 0.5$, $\sigma_2 \in [0.05,1]$. (b) Fixed $\sigma_1 = 1$, $\sigma_2 \in [0.05,1]$. (c) $\sigma_1 \in [0.05,1]$, fixed $\sigma_2 = 0.5$. (d)$\sigma_1 \in [0.05,1]$, fixed $\sigma_2 =1$. }
\label{fig:6}
\end{figure}

Furthermore,  inspired by \cite{menck2013basin} and \cite{serdukova2016stochastic},  we would like to employ  the stochastic basin of attraction (SBA) to quantify the basin stability in the escape region $D$.
 We denote the stochastic basin of the attractor $f_3$ in the region $D$ is the set $K( p^*) = \{(v,w) \in D| p(v,w) > p^*\}$, where $p^*$ indicates a high probability level (`threshold').
  The set $K(p^*)$  means the orbits starting from the points in the region $D$  reach to the region $E$ with high probability (we ignore the initial points whose solution have a `small' probability away from by $f_3$). The basin stability in region $D$ can be quantified in terms of its area, that is, $S_{FEP}$ represents the area of $K(p^*)$. More specifically, when the FEP of the solution orbit from the points in $D$ is greater than the threshold $p^*$, we record these points and then calculate the area of the region composed of these points as $S_{FEP}$. The normalized $S_{FEP}$ is \cite{CR2017FHN}
\begin{equation}\label{area}
  R_{FEP}=\frac{S_{FEP}(p^*)}{S_D},
\end{equation}
where $S_D$ is the area of region $D$.

Now we choose $p^*=0.8$ to compute $R_{FEP}$ as illustrated in Figure \ref{fig:5}. In Figure \ref{fig:5}, we let the noise intensity $\sigma_1=\sigma_2=\sigma$. From Figure \ref{fig:5}(a), we can see that for a fixed noise
intensity, the larger L\'evy motion index, the larger $R_{FEP}$. When $\alpha$ being 2 (corresponding to the case of
Brownian motion), the four curves all get their maximum and present a rapid growth. The situation shows that compared with the effect of L\'evy noise on the ML system, there are more orbits starting from the points in the low potential region to escape to the high potential region with high probability in the case of Brownian motion.
 It also indicates that Gaussian noise makes the ML system more accessible to the excited state. In Figure \ref{fig:5}(b), as can be seen, $R_{FEP}$ decreases with the increase of $\sigma$ for fixed $\alpha$ besides $\alpha$ being 2. When $\alpha$ being 2, as the noise intensity $\sigma$ increases, $R_{FEP}$ increases monotonously, reaches its maximum at $\sigma_b$, and then decreases as $\sigma$ continues to grow. As a whole, the lager L\'evy motion index and the smaller noise intensity are beneficial to more solution orbits with the point in region $D$ as the initial point escape to the target region $E$.

The pictorial representation of the aforementioned results
gives us an inspiration for the stochastic escape problem, i.e., if we expect the
ML system to escape from the resting state to the excited
state, then a smaller noise intensity and a larger L\'evy motion index $\alpha$
(smaller jump magnitude with higher frequency) should be
selected.

Next, we denote a new parameter named noise intensity ratio
\begin{equation}
r=\sigma_2 / \sigma_1.
\end{equation}
In order to better understand the influence of noise intensity and L\'evy motion index on FEP starting from the equilibrium point, we plot FEP depending on L\'evy motion index $\alpha$ and noise intensity ratio $r$, as shown in Figure \ref{fig:6}. In Figure \ref{fig:6}(a)-(b), we fix separately $\sigma_1=0.5$ and
$\sigma_1=1$, $\sigma_2$ belongs to interval $[0.05,1]$. It can be seen that the first escape probabilities all have lager value for larger L\'evy motion index $\alpha$ and smaller noise intensity ratio $r$.
And when $\sigma_1$ is fixed, the smaller $\sigma_2$ and the larger $\alpha$, the more likely  the ML system is to generate spikes. While in Figure \ref{fig:6}(c)-(d), we fix separately $\sigma_2=0.5$ and $\sigma_2=1$, $\sigma_1$ belongs to interval $[0.05,1]$. It can be seen that, for fixed noise intensity ratio
$r$, FEP is larger for the lager L\'evy motion index $\alpha$ in Figure \ref{fig:6}(c), but FEP is smaller for lager $\alpha$ and smaller $\sigma_1$ (lager noise intensity ratio $r$) with fixed $\sigma_2=1$ in Figure
\ref{fig:6}(d).  And when $\sigma_2$ is fixed, $\sigma_1$ has little effect on the state transition of the ML system.  In general, we infer that the noise intensity $\sigma_2$ has a greater influence on escape probability for the same $\alpha$.

\section{Mean first exit time}

\begin{figure}
\begin{minipage}[!h]{0.49 \textwidth}
\centerline{\includegraphics[width=7cm,height=3.8cm]{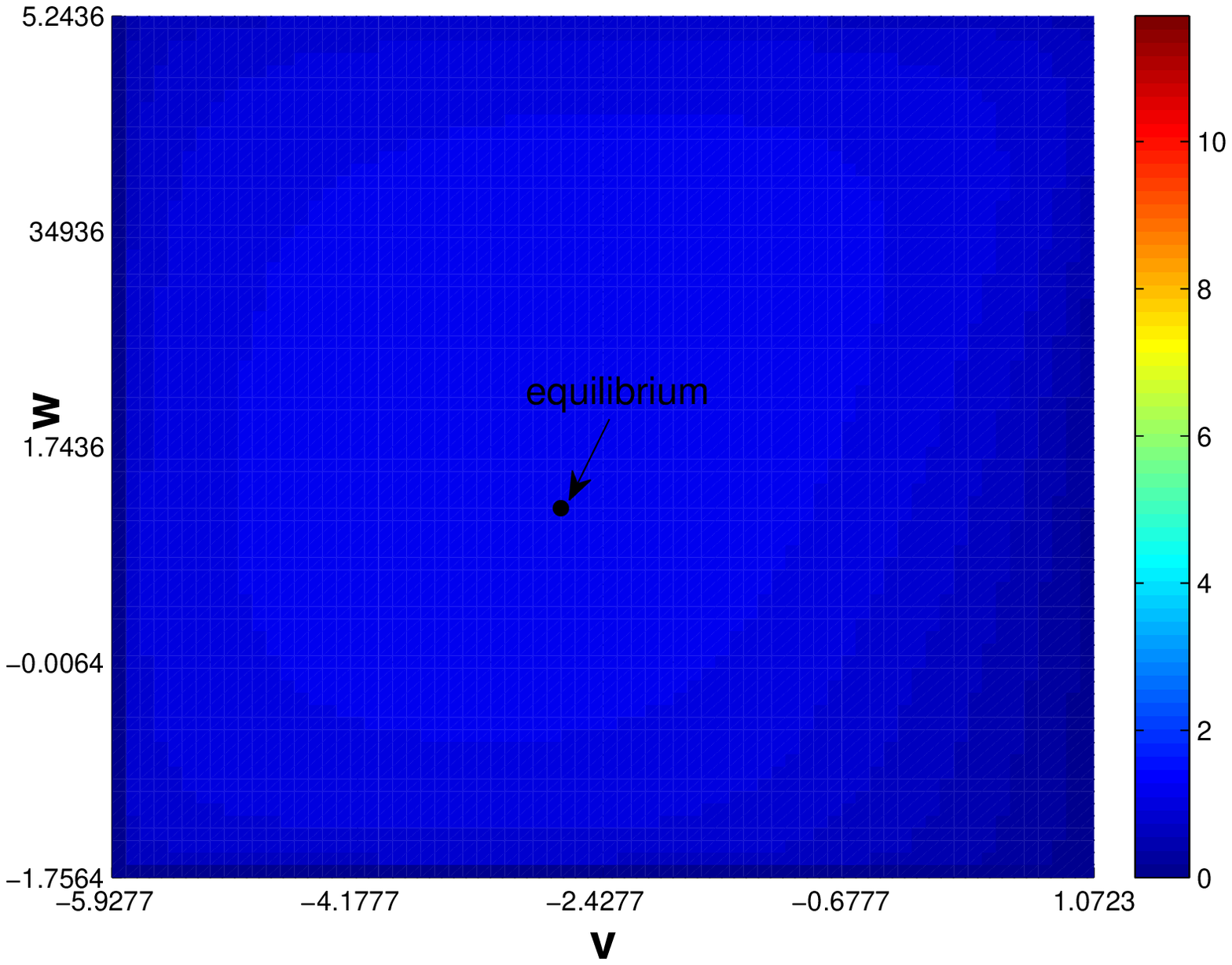}}
\centerline{(a) $\alpha$ = 0.5, $\sigma$ = 0.75.}
\end{minipage}
\hfill
\begin{minipage}[!h]{0.49 \textwidth}
\centerline{\includegraphics[width=7cm,height=3.8cm]{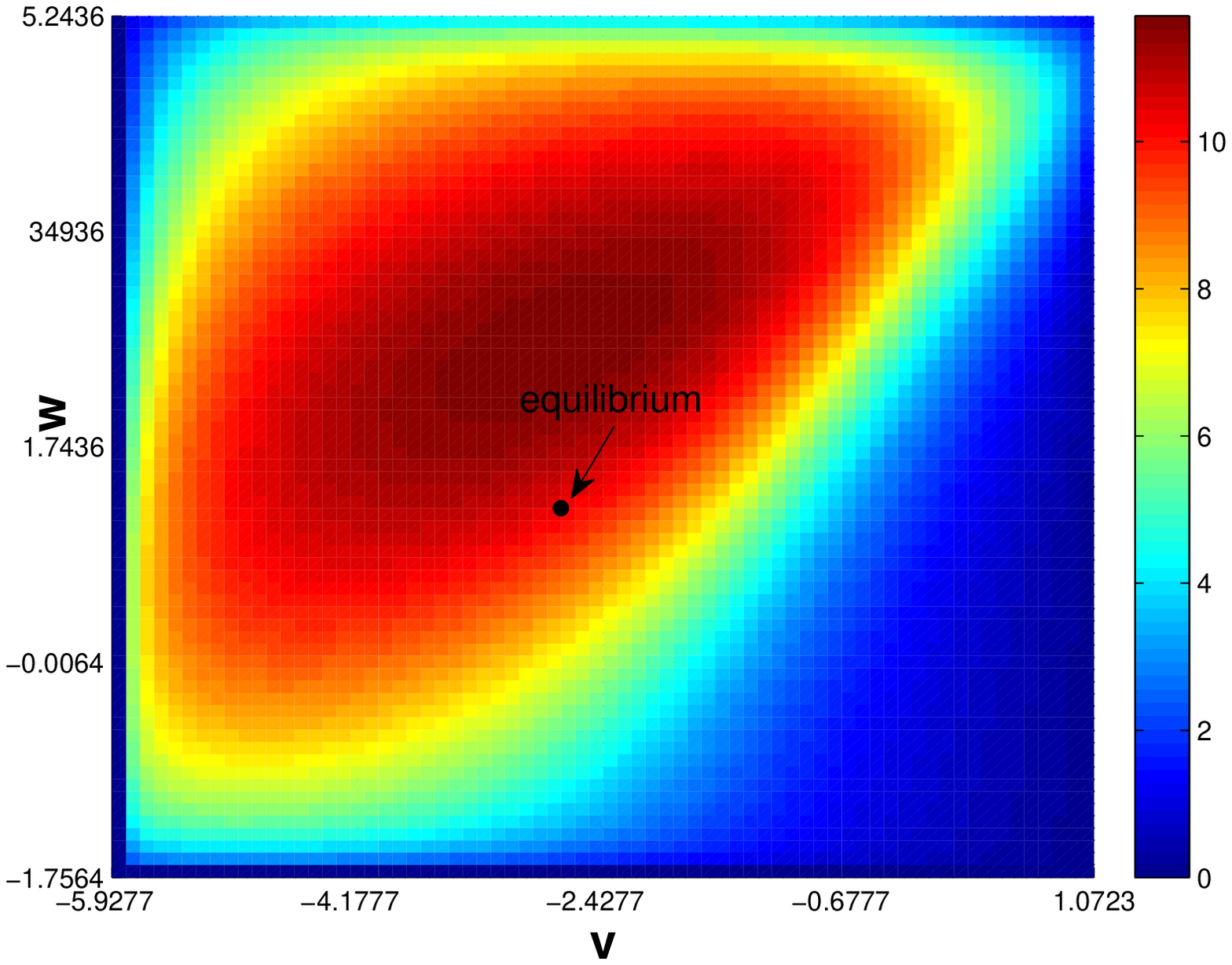}}
\centerline{(e) $\alpha$ = 1.25, $\sigma$ = 0.25.}
\end{minipage}

\begin{minipage}[!h]{0.49 \textwidth}
\centerline{\includegraphics[width=7cm,height=3.8cm]{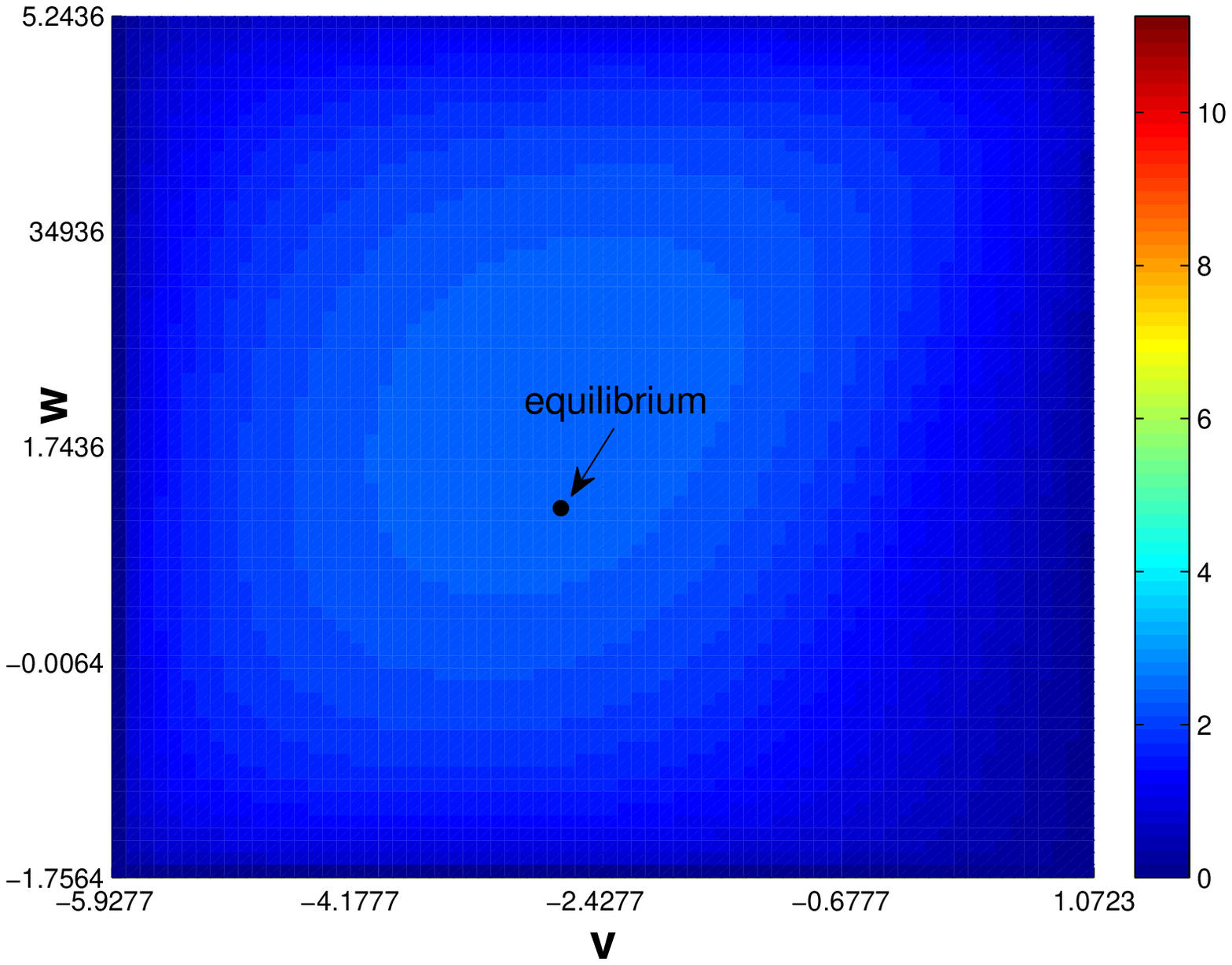}}
\centerline{(b) $\alpha$ = 1, $\sigma$ = 0.75.}
\end{minipage}
\hfill
\begin{minipage}[!h]{0.49 \textwidth}
\centerline{\includegraphics[width=7cm,height=3.8cm]{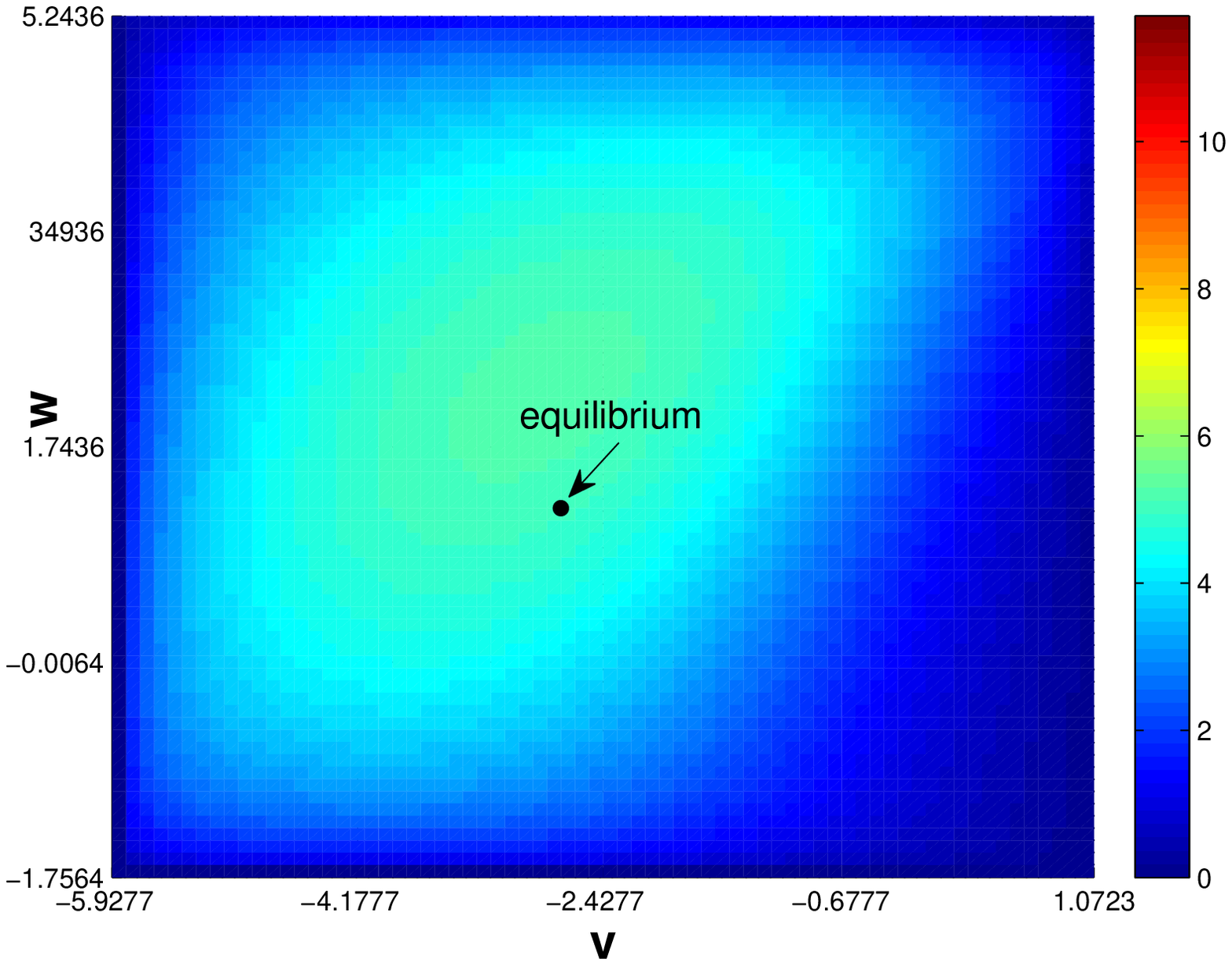}}
\centerline{(f) $\alpha$ = 1.25, $\sigma$ = 0.5.}
\end{minipage}

\begin{minipage}[!h]{0.49 \textwidth}
\centerline{\includegraphics[width=7cm,height=3.8cm]{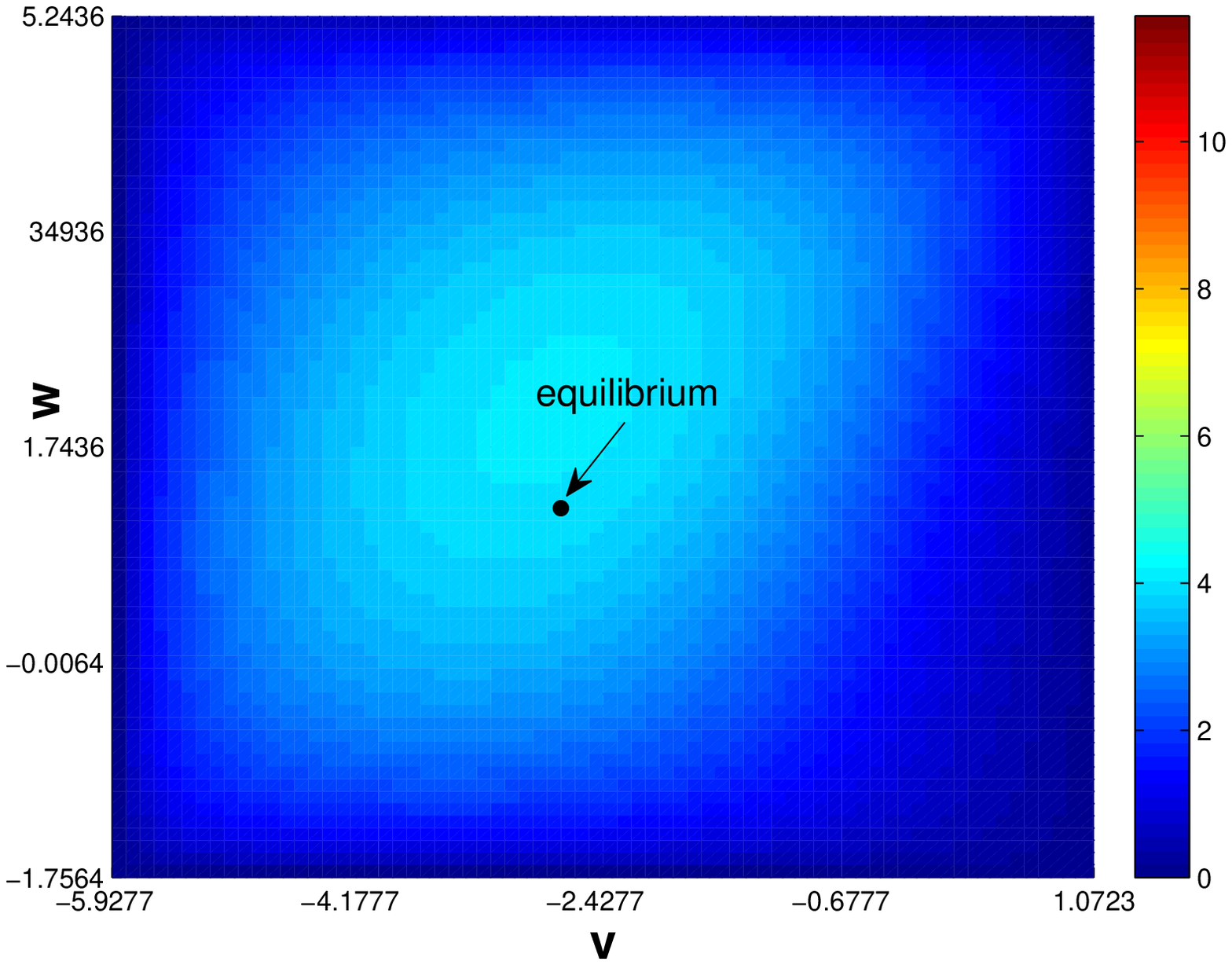}}
\centerline{(c) $\alpha$ = 1.5, $\sigma$ = 0.75.}
\end{minipage}
\hfill
\begin{minipage}[!h]{0.49 \textwidth}
\centerline{\includegraphics[width=7cm,height=3.8cm]{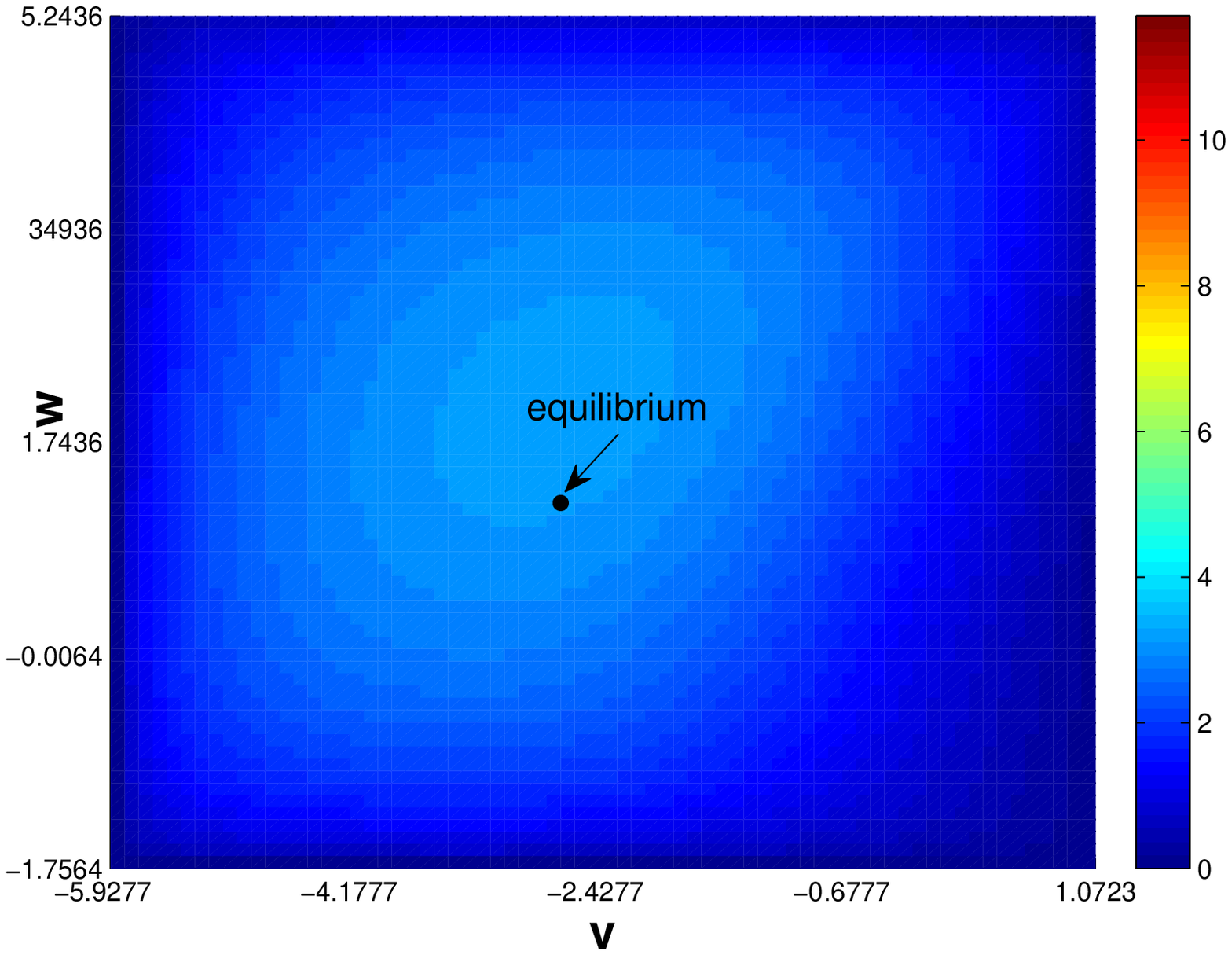}}
\centerline{(g) $\alpha$ = 1.25, $\sigma$ = 0.75.}
\end{minipage}

\begin{minipage}[!h]{0.49 \textwidth}
\centerline{\includegraphics[width=7cm,height=3.8cm]{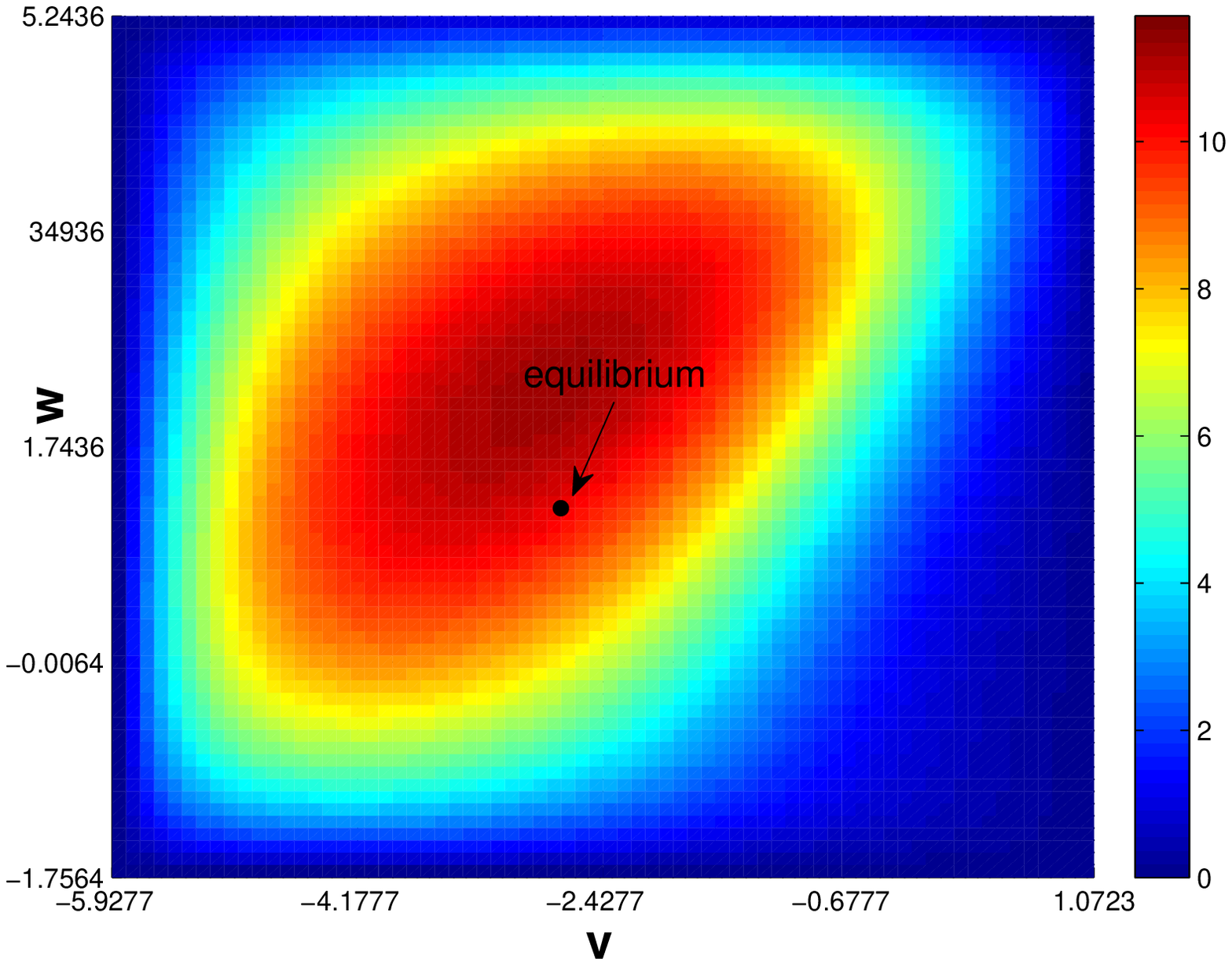}}
\centerline{(d) Brownain case, $\sigma$ = 0.75.}
\end{minipage}
\hfill
\begin{minipage}[!h]{0.49 \textwidth}
\centerline{\includegraphics[width=7cm,height=3.8cm]{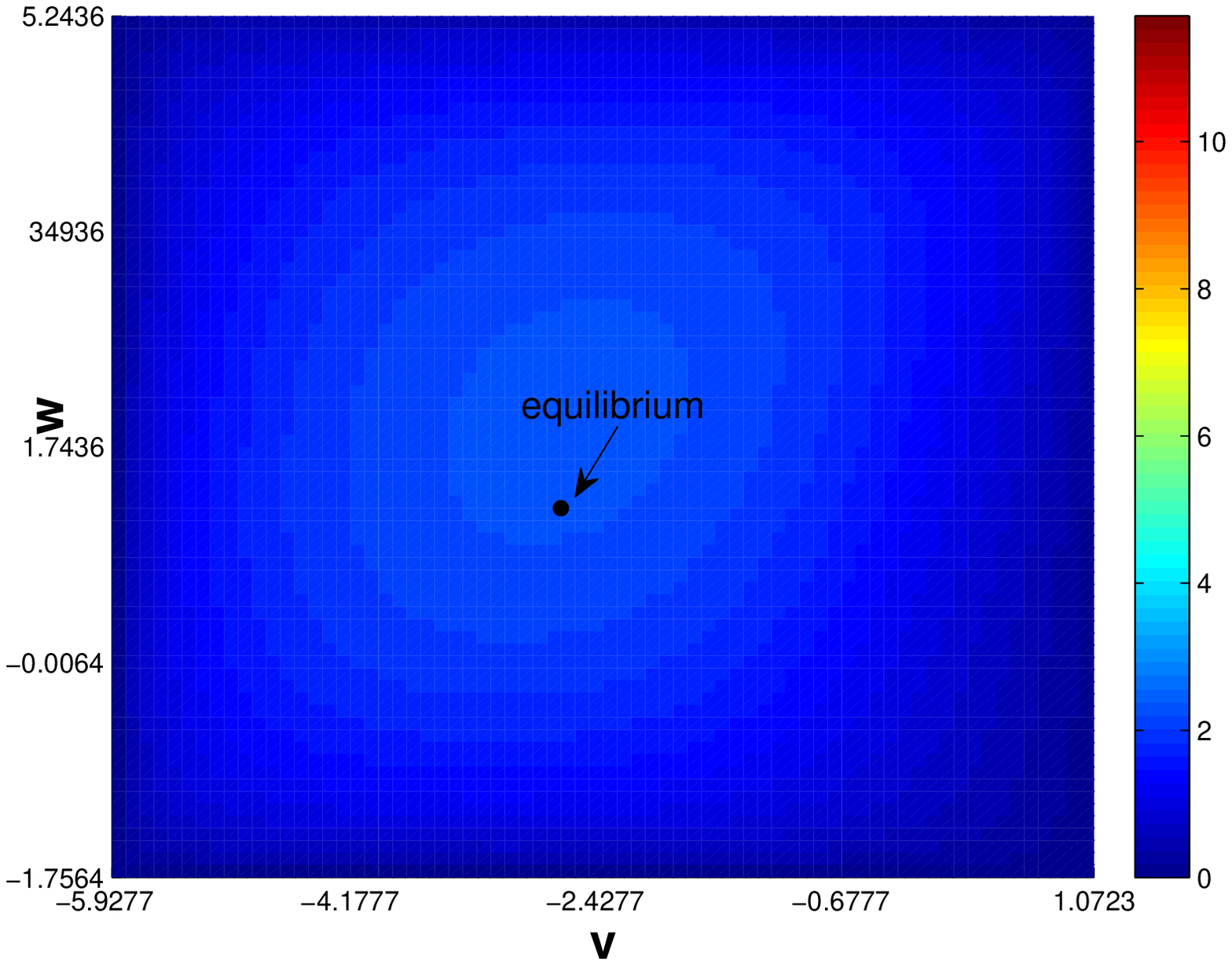}}
\centerline{(h) $\alpha$ = 1.25, $\sigma$ = 1.}
\end{minipage}

\vfill
\caption{MFET $u(v,w)$ from the escape region $D:(-5.9277,1.0723)\times(-1.7564,5.2436)$ to the region $D^{c}$. The colour map depends on L\'evy motion index $\alpha$ and noise intensity $\sigma$ ($\sigma_1 = \sigma_2 = \sigma$). (a)-(d) correspond to fixed noise intensity $\sigma=0.75$ and different L\'evy motion index $\alpha=0.5,~1,~1.5,~2$. (e)-(h) correspond to fixed $\alpha=1.25$ and different $\sigma=0.25,~0.5,~0.75,~1$. All figures are unified into an identical color map with the same scale, red marking 11.7084 and blue making 0.}
\label{fig:7}
\end{figure}

\begin{figure}[!ht]
\begin{minipage}{0.49 \textwidth}
\centerline{\includegraphics[width=9cm,height=6cm]{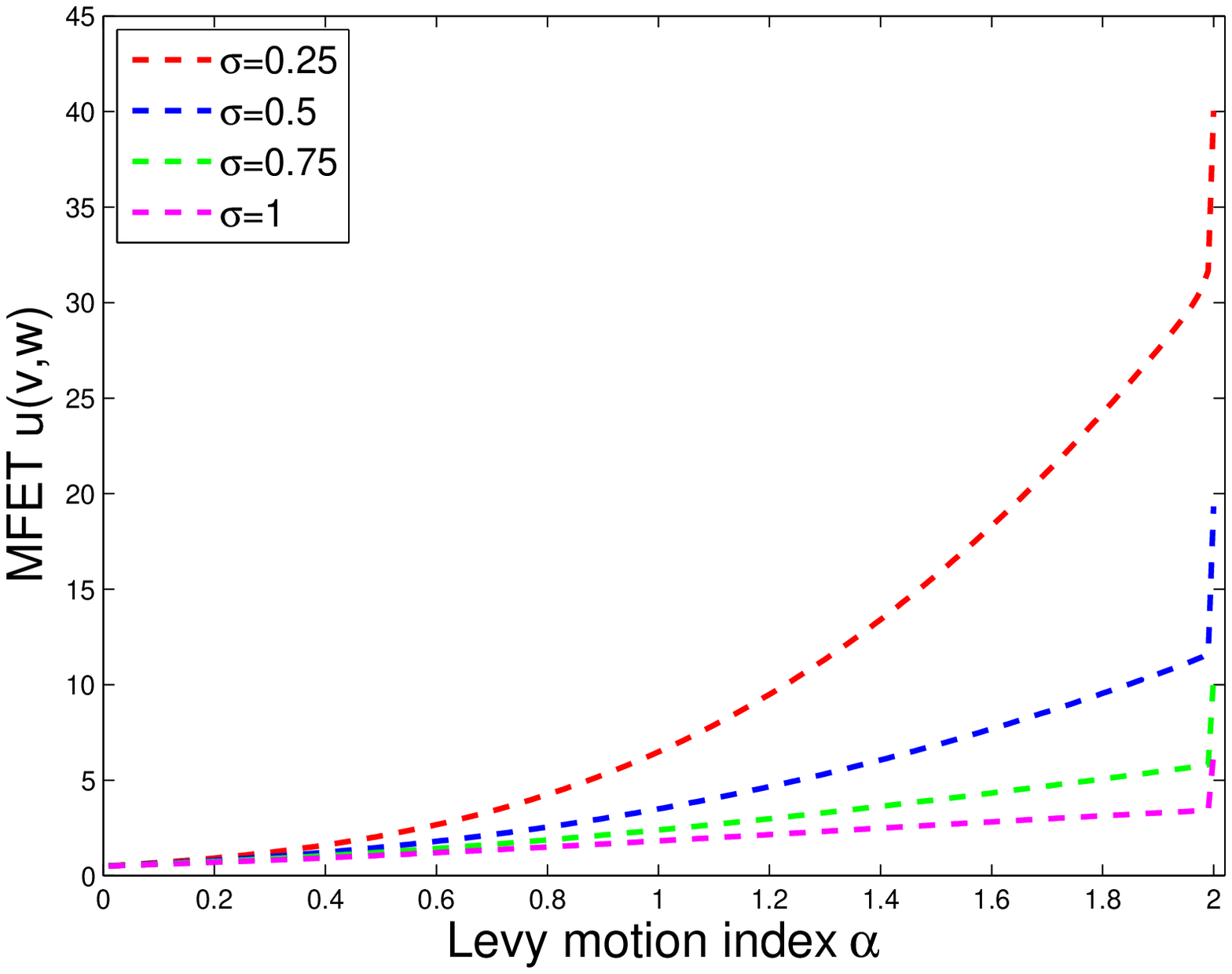}}
\centerline{(a)}
\end{minipage}
\hfill
\begin{minipage}{0.49 \textwidth}
\centerline{\includegraphics[width=9cm,height=6cm]{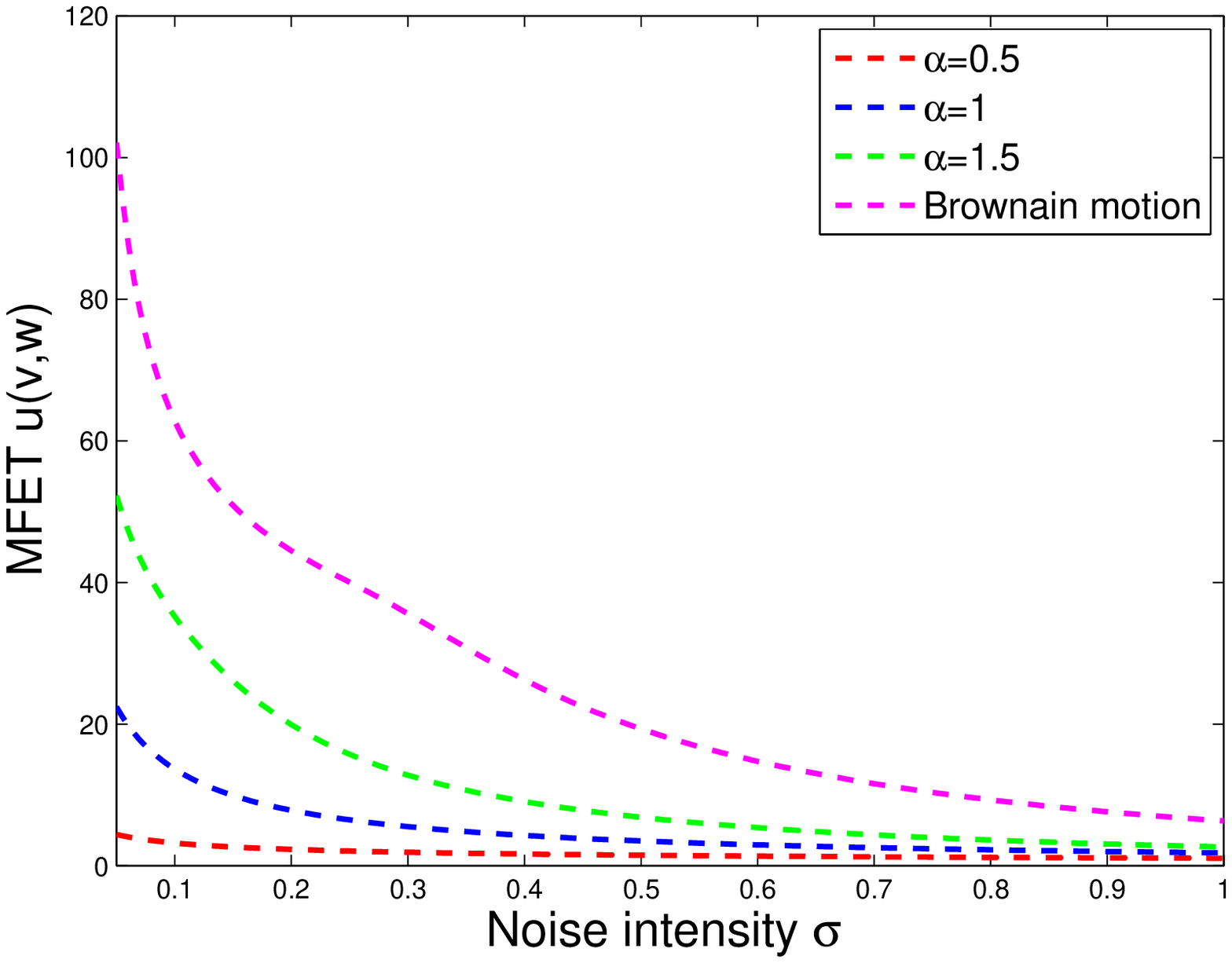}}
\centerline{(b)}
\end{minipage}

\vfill
\caption{MFET $u(v, w)$ of the equilibrium $s^*$= (-2.7277, 1.2436). Here $\sigma_1=\sigma_2=\sigma$. (a) Effect of L\'evy motion index $\alpha$ on MFET at $s^*$ with different noise intensity $\sigma$ (red: $\sigma$ = 0.25, blue: $\sigma$ = 0.5, green: $\sigma$ = 0.75,
pink: $\sigma$ = 1). (b) Effect of noise intensity $\sigma$ on MFET at $s^*$ with different L\'evy motion index $\alpha$ (red: $\alpha$ = 0.5, blue:
$\alpha$ = 1, green: $\alpha$ = 1.5, pink: $\alpha$ being 2 indicates to case of Brownain motion).}
\label{fig:8}
\end{figure}

In this section, we use another quantity: the mean first exit time  to examine the effects of L\'evy motion index and noise intensity on the behavior of the escape problem. Back to the general two-dimensional stochastic system (\ref{general}) in the previous section.
The first exit time for a solution orbit starting at $(v,w)$  in a region $D$ is defined as
\begin{equation}
  \tau(\omega,(v,w)) = \inf\{t \geq 0: (v_0,w_0) = (v,w), (v_t, w_t)  \in   D^c\},\notag
\end{equation}
where $ D^{c}$ is the complement set of the bounded region $D$ in $\mathbb{R}^{2}$.  The mean first exit time (MFET) is denoted as
\begin{equation}
  u(v,w):= \mathbb{E}\tau(\omega, (v, w)) \geq 0. \notag
\end{equation}

It is the `average' residence time of a solution orbit initially at $(v,w)$ inside region $D$ before escaping to another region. The difference between a Gaussian and a non-Gaussian process is that the orbit typically hits the boundary of region $D$ for the Brownian motion case, while it may jumps outside of region $D$ for L\'evy motion case. The MFET $u(v,w)$ satisfies the following nonlocal integral-differential equation with an exterior boundary condition \cite{duan2015introduction}:
\begin{eqnarray}\label{MFET}
  &Au(v,w)= -1, \qquad &(v,w) \in D,\notag \\
  &u(v,w)= 0,  \qquad &(v,w) \in D^{c}.
\end{eqnarray}
Here the generator $A$ is defined in equation \eqref{GNR}.  The existence and uniqueness of solution to the equation \eqref{MFET} satisfied by the mean first exit time of the stochastic ML system \eqref{SV} can be proved according to the third section in \cite{RO2016} and Theorem 3.2 in \cite{AA2016}. The equation \eqref{MFET}
can be solved by an effective numerical scheme given in the Appendix.

MFET can be used as a tool to measure the stability of a system: The longer the MFET, the more stable the resting state. For ML system, we choose the low potential region $D$ enclosing the resting state $s^*$ in phase plane to calculate the MFET. We plot MFET $u(v,w)$ from escape region $D:(-5.9277,1.0723)\times(-1.7564,5.2436)$ to the region $D^{c}$ as shown in Figure \ref{fig:7}. In Figure \ref{fig:7}(a)-(d), noise intensity $\sigma$ is fixed as $\sigma=0.75$, the MFET of non-exciting region gradually increases with the increase of L\'evy motion index $\alpha$ ($\alpha$= 0.5,~1,~1.5,~2). In Figure \ref{fig:7}(e)-(h), as $\sigma$ increases, the MFET of non-exciting region is getting less for fixed $\alpha=1.25$. In fact, the solution orbit starting at the equilibrium point will stay there forever without noise. Now for noisy situations, the solution orbits may stay in region $D$ for a finite time and then escape to another region.
So we can employ the MFET to characterize the relative stability of the solution orbits starting at the points in region $D$. From Figure \ref{fig:7}(a)-(h), we observe that larger L\'evy motion index $\alpha$ and smaller noise intensity $\sigma$ are beneficial to the stability of the region $D$.

\begin{figure}[!ht]
\begin{minipage}{0.49 \textwidth}
\centerline{\includegraphics[width=9cm,height=6cm]{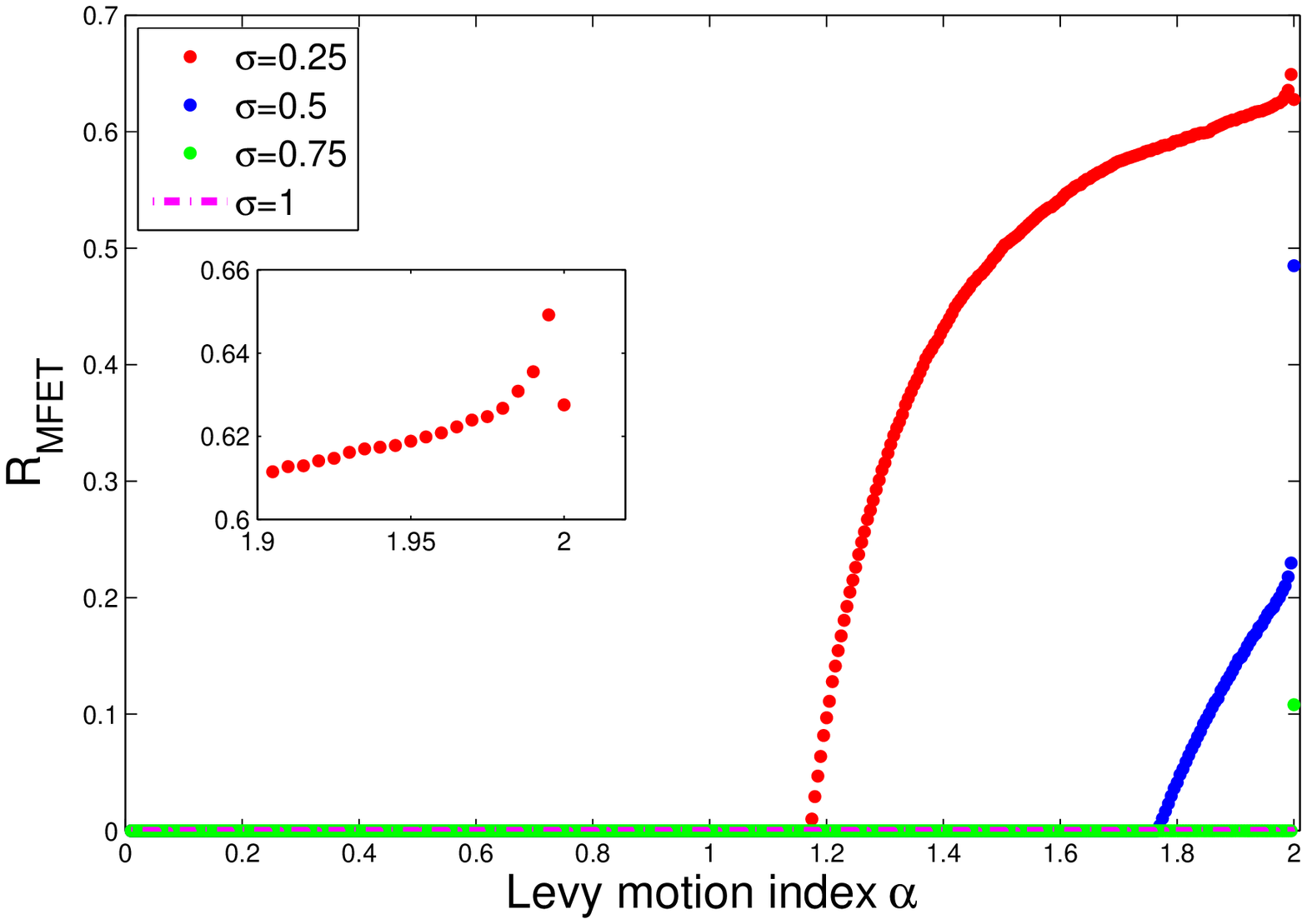}}
\centerline{(a)}
\end{minipage}
\hfill
\begin{minipage}{0.49 \textwidth}
\centerline{\includegraphics[width=9cm,height=6cm]{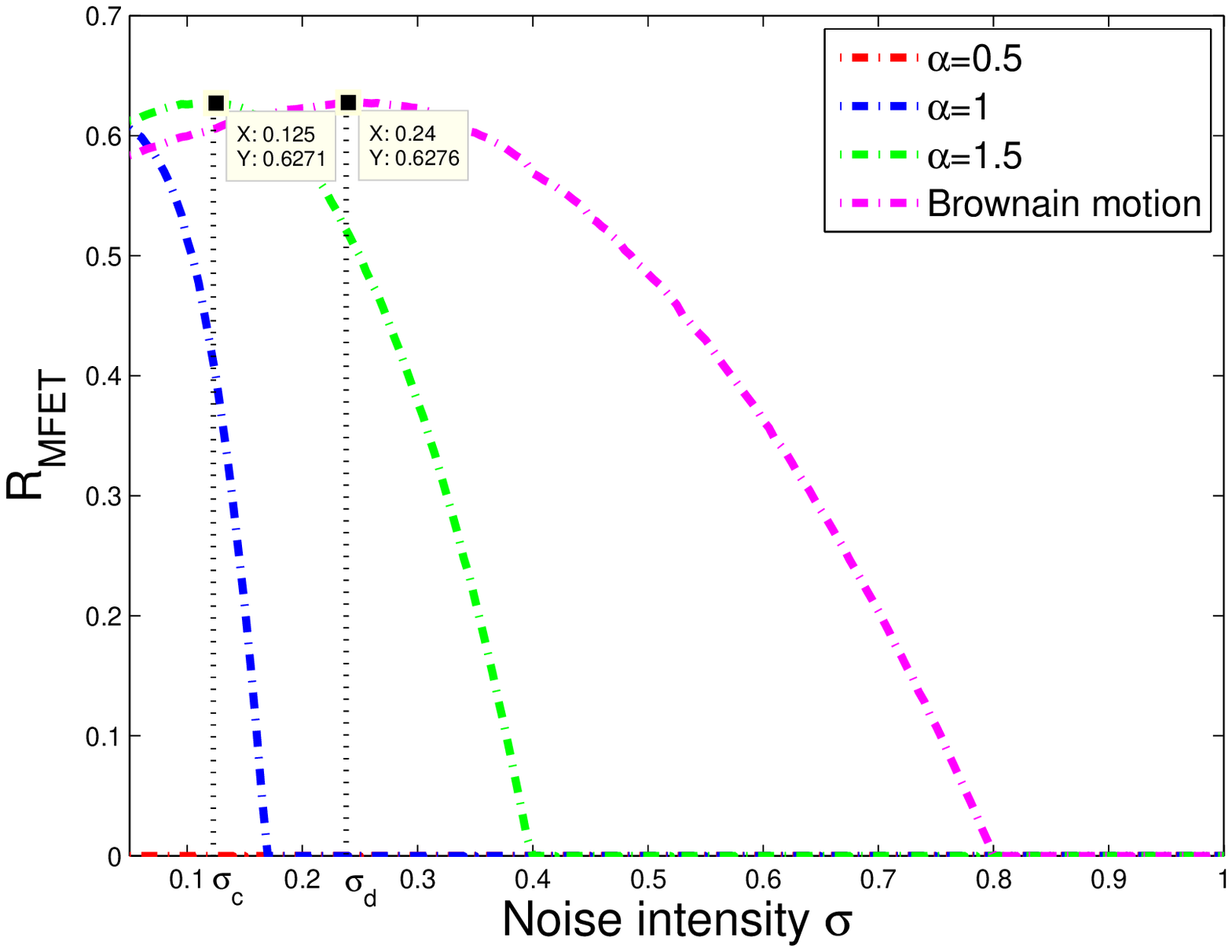}}
\centerline{(b)}
\end{minipage}

\vfill
\caption{Effect of noise intensity $\sigma$ and L{\'e}vy motion index $\alpha$ on the $R_{MFET}$. The
threshold $u^*$ to compute $R_{MFET}$ is chosen as 10. Here we choose $\sigma_1=\sigma_2=\sigma$. (a) $R_{MFET}$ against L\'evy motion
index $\alpha$ for various noise intensity $\sigma$ (red: $\sigma$ = 0.25, blue $\sigma$ = 0.5, green $\sigma$ = 0.75,
pink $\sigma$ = 1). (b) $R_{MFET}$ against noise intensity $\sigma$ for various L{\'e}vy motion index $\alpha$ (red: $\alpha$ = 0.5, blue:
$\alpha$ = 1, green: $\alpha$ = 1.5, pink: $\alpha$ indicates to case of Brownain motion).}
\label{fig:9}
\end{figure}

Similarly, Figure \ref{fig:8} depicts the change of MFET for the resting state in some cases. In Figure \ref{fig:8}, we denote the L\'evy motion index $\alpha \in (0,2)$ (the case where $\alpha$ being 2 is replaced by the Brownian motion) and the noise intensity $\sigma \in (0, 1]$ ($\sigma_1=\sigma_2=\sigma$).
It can be seen from the Figure \ref{fig:8}(a), the larger $\alpha$, the greater MFET, and this is more noticeable for smaller noise intensity. And when the noise intensity $\sigma$ is larger ($\sigma=1$), the L\'evy motion $\alpha$ is less effect on the MFET for the resting state. No matter what the noise intensity is ($\sigma=0.25,~0.5,~0.75,~1$), MFET for the resting state reaches the maximum at $\alpha$ being 2. The MFET for the Brownian motion is longer than for the L\'evy motion, which may be due to the fact that the trajectory of Brownian motion is continuous and cannot jump. This also indicates the resting state is more stable under Brownian motion than under L\'evy motion.

 Indeed,  MFET is solution to the \eqref{MFET}, with $A$ in \eqref{GNR} for L\'evy case and  \eqref{GNRbm} for Brownain case. In Figures 8(a), the numerical simulation indicates discontinuity as $\alpha$  approaching 2. Theoretically, although the characteristic function of L\'evy motion $L_{t}^{\alpha}$ satisfies $e^{-t|\xi|^{\alpha}} \rightarrow e^{-t|\xi|^{2}}$ as $\alpha \rightarrow 2$, the MFET  for L\'evy case may not approach that for Brownian case as $\alpha \rightarrow 2$; see Theorem 2.2 in \cite{Imkeller2009}. This remark also applies to FEP.

In addition, Figure \ref{fig:8}(b) further depicts the effect of noise intensity $\sigma$ on MFET for the resting state with different L\'evy motion index $\alpha$. As can be seen, MFET decreases with the increasing of $\sigma$ for fixed $\alpha$ and the change curve of MFET is not very obvious when $\alpha$ is small, such as $\alpha=0.5$. On the contrary, when $\alpha$ is lager, the change curve of MFET goes down very fast at first and tends to level off with the increasing of $\sigma$.  The noise intensity has little effect on the MFET for fixed smaller L\'evy motion index $\alpha$ ($\alpha=0.5$).
In general, the smaller the noise intensity $\sigma$ and the larger the L\'evy motion $\alpha$, the longer the time for the orbit to escape the region $D$. Therefore, if we expect
the resting state to be more stable, then a smaller noise intensity $\sigma$ and a larger $\alpha$ (smaller jump magnitude with higher frequency) should be responsible.

We now define another stability concept via to the stochastic basin of attraction  in Section 3.
Let $M(u^*) = \{(v,w) \in D|u(v,w) > u^* \}$, i.e., the solution orbit starting
from region $D$ and remaining there for a finite time (remarked by a threshold $u^*$). Then we also
quantify the basin stability in region $D$ based on its area \cite{CR2017FHN}:
\begin{equation}\label{area1}
  R_{MFET}=\frac{S_{MFET}(u^*)}{S_D},
\end{equation}
where $S_{MFET}$ is the area of $ M(u^*)$, and $S_D$ is the area of $D$. $R_{MFET}$ is the normalization of $S_{MFET}$.

\begin{figure}[!ht]
\begin{minipage}{0.49 \textwidth}
\centerline{\includegraphics[width=9cm,height=6cm]{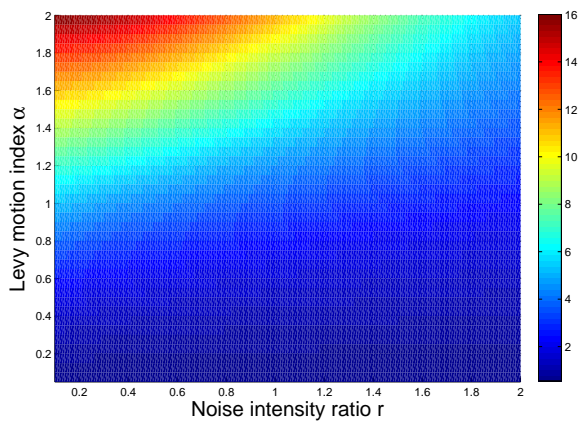}}
\centerline{(a) $r=\sigma_2/\sigma_1$, fixed $\sigma_1=0.5$.}
\end{minipage}
\hfill
\begin{minipage}{0.49 \textwidth}
\centerline{\includegraphics[width=9cm,height=6cm]{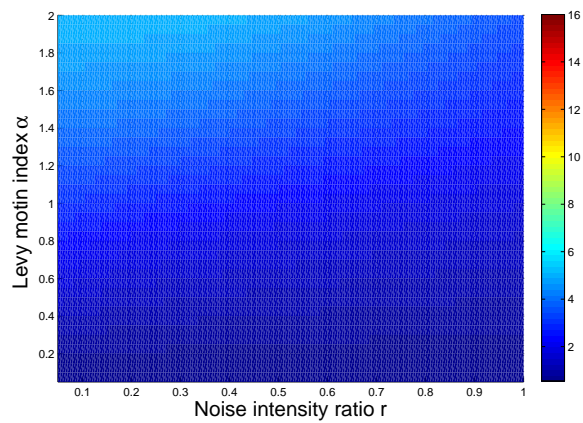}}
\centerline{(b) $r=\sigma_2/\sigma_1$, fixed $\sigma_1=1$.}
\end{minipage}

\begin{minipage}{0.49 \textwidth}
\centerline{\includegraphics[width=9cm,height=6cm]{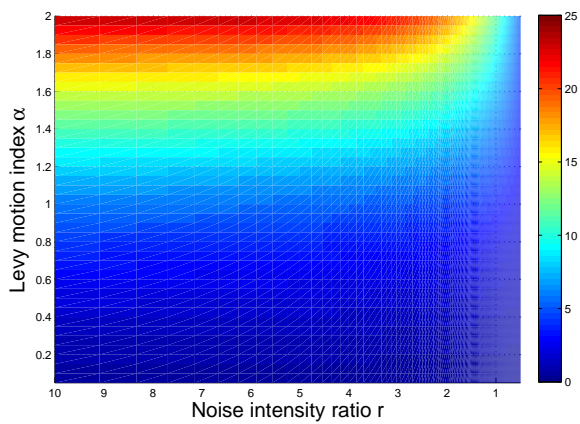}}
\centerline{(c) $r=\sigma_2/\sigma_1$, fixed $\sigma_2=0.5$.}
\end{minipage}
\hfill
\begin{minipage}{0.49 \textwidth}
\centerline{\includegraphics[width=9cm,height=6cm]{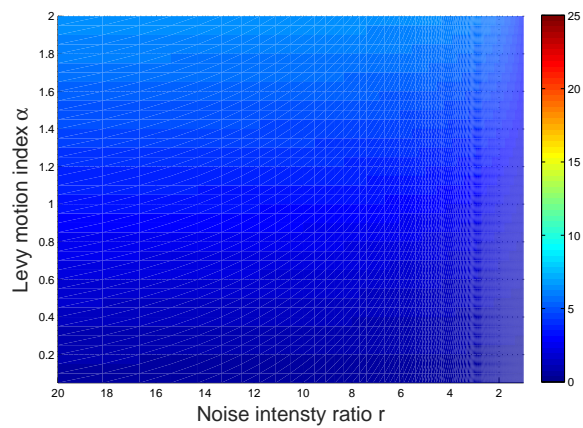}}
\centerline{(d) $r=\sigma_2/\sigma_1$, fixed $\sigma_2=1$.}
\end{minipage}

\vfill
\caption{MFET $u(v,w)$ of the equilibrium $s^*$ = (-2.7277, 1.2436), the color graphs depend on  noise intensity ratio $r$ and L\'evy motion index $\alpha$. (a)-(b) $\sigma_2 \in [0.05,1]$, fixed $\sigma_1=0.5,~1$, respectively. The color bar of (a) and (b) is taken the same scale. (c)-(d) $\sigma_1 \in [0.05,1]$, fixed $ \sigma_2=0.5,~1$, respectively. The color bar of (c) and (d) is taken the same scale.}
\label{fig:10}
\end{figure}

In the following, as an example, we make $u^*=10$ to calculate the $R_{MFET}$, the results are shown in Figure \ref{fig:9}. Figure \ref{fig:9}(a) depicts $R_{MFET}$ against the L\'evy motion index $\alpha$ for various values of $\sigma$ ($\sigma_1=\sigma_2=\sigma$). It can be seen that for fixed $\sigma=0.25$ and $\sigma=0.5$, the $R_{MFET}$ remains 0 in the beginning, then increases with increasing $\alpha$ besides $\alpha$ being 2 (corresponding to the case of Brownian motion). While for the other two curves of
$R_{MFET}$ with $\sigma=0.75$ and $\sigma=1$, the $R_{MFET} = 0$, which means that for higher noise intensity the solution orbit starting at a point in region $D$ gets out of region $D$ quickly. This also indicates the region $D$ is less stable under the higer noise intensity. Special cases occur when $\alpha$ being 2: for $\sigma=0.5$ and $\sigma=0.75$, the values of $R_{MFET}$ both have a jumping growth; for $\sigma=0.25$,
$R_{MFET}$ suddenly decreases and the value of $R_{MFET}$ has no change for $\sigma=1$. Figure \ref{fig:9}(b) depicts $R_{MFET}$ against noise intensity $\sigma$ for various values of $\alpha$. It shows that for $\alpha=1.5$ and $\alpha$ being 2, the curve of $R_{MFET}$ transits from increasing to decreasing at $\sigma_c$ and $\sigma_d$, respectively, then goes down to zero as further increasing $\sigma$. While when $\alpha=1$, $R_{MFET}$ gradually decreases from positive to zero, then keep zero with the further increase of $\sigma$. At the same time, $R_{MFET}=0$ when $\alpha=0.5$, which indicates that a solution orbit
starting at $(v,w) \in D $ escapes quickly for smaller $\alpha$, that is, the region $D$ is less stable for the smaller L\'evy motion index. Here we choose the initial value $\sigma=0.05$ and the same partition $\Delta \sigma=0.005$. Moreover, we can find that the curve with $\alpha = 2$ and $\alpha = 1.5$, $\alpha=1$ has a intersection point, respectively. For fixed noise intensity, larger $\alpha$ has smaller $R_{MFET}$ before the intersection point and it's the opposite after the intersection point.

As in the section 3, we also plot MFET for the resting state depending on noise intensity ratio $r$ and L\'evy motion index $\alpha$ as shown in the Figure \ref{fig:10}. In Figure \ref{fig:10}(a)-(b), we respectively fix
$\sigma_1=0.5$ and $\sigma_1=1$ as well $\sigma_2$ belongs to interval $[0.05,1]$. It can be seen that MFET has lager value for lager L\'evy motion index $\alpha$ and smaller noise intensity ratio $r$. Meanwhile, when fix
$\alpha<1$, whether in Figure \ref{fig:10}(a) or Figure \ref{fig:10}(b), the noise intensity ratio $r$  has a
small effect on MFET of the equilibrium point; and when fix $\alpha>1$, the change of MFET is obvious.
This just verifies that the jump property of L\'evy motion. While in Figure \ref{fig:10}(c)-(d), we respectively fix $\sigma_2=0.5$ and $\sigma_2=1$ as well $\sigma_1$ belongs to interval $[0.05,1]$.  We observe that the lager $\alpha$, the longer MFET and the noise intensity $\sigma_1$ has less effect on the system for the fxed $\sigma_2$ ($\sigma_2=1$). Comparing Figure \ref{fig:10}(c) with Figure \ref{fig:10}(d) , when fix $\alpha>1$, for the same $\sigma_1$ (the $r$ in (d) is 2 times that in (c)), the MFET with $\sigma_2=0.5$ obviously larger than the MFET with $\sigma_2=1$, which indicates the resting state is relatively stable for smaller $\sigma_2$. In general, the larger the L\'evy motion index and the smaller the noise intensity, the more stable the system.

\section{Conclusion}

In summary, we focus on the escape problem driven by a symmetric $\alpha$-stable L\'evy noise (non-Gaussian noise) in Morris-Lecar (ML) model. We have provided a method to quantify  the dynamics of escape from the resting state of the system by means of two deterministic indices: the first escape probability (FEP) and the mean first exit time (MFET). To be specific, we have used the method to describe the state transition from the resting state to the excited state. To simulate the firing behavior of neurons and calculate FEP, we have chosen the appropriate escape region containing the equilibrium point and the target region. Meanwhile, we have depicted the stability of the escape region in terms of  MFET.

Through numerical simulation and analysis, we have found that the noise intensities $\sigma_1$, $\sigma_2$ and the jump size of $\alpha$-stable L\'evy motion have significant and delicate influences on the FEP and MFET.
We have also discovered  that for smaller jumps of the L\'evy motion and relatively smaller noise
intensity, FEP is larger, which means that they are conducive to the production of spikes.
 However, higher noise intensity and larger jumps of the L\'evy motion shortens the MFET,  which means the escape region is less stable for higher noise intensity and larger jumps of the L\'evy motion. Moreover, the Brownian motion (Gaussian noise) has been also considered and compared with the L\'evy motion case.
Compared with L\'evy motion, FEP is larger and MEFT is longer under the influence of Brownian motion. This means that the perturbed system is more likely to switch from the resting state to the excited state, and the resting state is relatively more stable in the Brownian case. Meanwhile, we have explored the stochastic basin of attraction for the stability of a region.
By calculating the area of high FEP and long MFET for various L\'evy motion index $\alpha$ and noise intensity $\sigma$, we have revealed that larger $\alpha$ and smaller noise intensity  $\sigma$ are beneficial for the stability of the region $D$.

By calculating the impact of  the noise intensity ratio $r$ and L\'evy motion index $\alpha$ on FEP of the equilibrium point,
we have revealed  that $\sigma_2$  (intensity for ion channel noise) has more pronounced  influence on the system than $\sigma_1$ (intensity for current fluctuations),  for fixed L\'evy motion index $\alpha$. To be specific, for the state transition probability of the system, ion channel noise has a greater impact than current noise. The smaller the $\sigma_2$ and the larger L\'evy motion index $\alpha$, the more likely  the stochastic ML system is to generate a pulse. For the MFET at the equilibrium point, larger L\'evy motion index $\alpha$ and smaller noise intensities  $\sigma_1$, $\sigma_2$, the more stable the resting state in the ML system.

We have applied the knowledge of stochastic
dynamics to explore the phenomenon of noise-induced escape in a neural system. This work provides some  mathematical understanding about the impact of non-Gaussian, heavy-tailed, burst-like fluctuations on excitable systems such as the Morris-Lecar system.

\section*{Acknowledgements}
We would like to thank  Xiaoli Chen, Wei Wei and Yongge Li for helpful discussions. This work was partly supported by the National
Science Foundation Grant (NSF) No. 1620449, and the National Natural Science Foundation of China (NSFC) Grant Nos. 11531006 and 11771449.

\section*{Appendix: Numerical simulation}
 \label{SI}

\renewcommand{\theequation}{A.\arabic{equation}}
\renewcommand{\thefigure}{A.\arabic{figure}}
\setcounter{equation}{0}
\setcounter{figure}{0}

We use an efficient numerical finite difference scheme  \cite{gao2014mean} to compute the FEP (equation (\ref{FEP})) and MFET (equation (\ref{MFET})). This method is revised for our model in $(v,w)\in D=(a,b)\times(c,d)$, $E=[a',b']\times[c,d]$ by a scalar conversion $v=\frac{b-a}{2}s+\frac{a+b}{2}$, and $w=\frac{d-c}{2}k+\frac{c+d}{2}$ for $s\in(-1,1),k\in(-1,1)$. If we let $m(s,k)=p(\frac{b-a}{2}s+\frac{a+b}{2},\frac{d-c}{2}k+\frac{c+d}{2})$, equation (\ref{GNR}) is discretized as follows:
\begin{align}
Ap(v,w)=&\frac{2}{b-a}f_{1}(\frac{b-a}{2}s+\frac{a+b}{2},\frac{d-c}{2}k+\frac{c+d}{2})m_{s}\nonumber\\
&+\frac{2}{d-c}f_{2}(\frac{b-a}{2}s+\frac{a+b}{2},\frac{d-c}{2}k+\frac{c+d}{2})m_{k}\nonumber\\
&-\frac{\sigma_{1}^{\alpha}C_{\alpha}}{\alpha}(\frac{b-a}{2})^{\alpha}\left[\frac{1}{(1+s)^{\alpha}}
+\frac{1}{(1-s)^{\alpha}}\right]m(s,k) \nonumber\\ &+\sigma_{1}^{\alpha}C_{\alpha}(\frac{b-a}{2})^{\alpha}\int_{-1-s}^{1-s}\frac{m(s+s',k)-m(s,k)}{\mid s'\mid^{1+\alpha}}ds'  \nonumber\\
&-\frac{\sigma_{2}^{\alpha}C_{\alpha}}{\alpha}(\frac{d-c}{2})^{\alpha}\left[\frac{1}{(1+k)^{\alpha}}
+\frac{1}{(1-k)^{\alpha}}\right]m(s,k) \nonumber\\ &+\sigma_{2}^{\alpha}C_{\alpha}(\frac{d-c}{2})^{\alpha}\int_{-1-k}^{1-k}\frac{m(s,k+k')-v(s,k)}{\mid k'\mid^{1+\alpha}}dk' \nonumber\\
=& \psi(s,k), \nonumber\\
\end{align}
where the integral in this equation is taken as the Cauchy principle value integral, and $$\psi(s,k)=\frac{\sigma_{1}^{\alpha}}C_{\alpha}{\alpha}\left[\frac{1}{(b'-\frac{b-a}{2}s-\frac{a+b}{2})^{\alpha}}
-\frac{1}{(a'-\frac{b-a}{2}s-\frac{a+b}{2})^{\alpha}}\right]$$
in the case of FEP (equation (\ref{FEP})), or $\psi(s,k)=-1$ in the case of MFET (equation (\ref{MFET})).




\end{document}